\documentclass[a4paper]{amsart}
\usepackage{amssymb} 
\usepackage[T1]{fontenc}
\usepackage[latin1]{inputenc}
\usepackage{amsfonts}
\usepackage{amsxtra}
\usepackage{ae}
\usepackage[all]{xy}
\usepackage{enumerate}

\include{diagram}

\newcommand*{\ket}{\rangle}
\newcommand*{\bra}{\langle}
\newcommand*{\ad}{\mathsf{ad}}

\newcommand*{\A}{\mathcal{A}}
\newcommand*{\B}{\mathcal{B}}

\newcommand*{\E}{\mathcal{E}}

\renewcommand*{\H}{\mathcal{H}}

\newcommand*{\K}{\mathcal{K}}
\renewcommand*{\L}{\mathcal{L}}
\newcommand*{\M}{\mathcal{M}}
\newcommand*{\N}{\mathcal{N}}
\renewcommand*{\P}{\mathcal{P}}

\renewcommand*{\S}{\mathcal{S}}

\newcommand*{\X}{\mathcal{X}}
\newcommand*{\CC}{\mathcal{CC}}
\newcommand*{\TC}{\mathcal{TC}}
\newcommand*{\CI}{\mathcal{CI}}
\newcommand*{\TI}{\mathcal{TI}}
\newcommand*{\cotimes}{\hat{\otimes}}
\newcommand*{\twisted}{\boxtimes}

\newcommand*{\DD}{\mathsf{D}}

\newcommand*{\Corep}{\mathsf{Corep}}
\newcommand*{\Irr}{\mathsf{Irr}}
\renewcommand*{\max}{\mathsf{f}}
\newcommand*{\red}{\mathsf{r}}

\newcommand*{\cop}{\mathsf{cop}}

\newcommand*{\CH}{\mathbb{C}}
\newcommand*{\FH}{\mathbb{F}}
\newcommand*{\HH}{\mathbb{H}}
\newcommand*{\KH}{\mathbb{K}}
\newcommand*{\LH}{\mathbb{L}}

\DeclareMathOperator{\coker}{coker}
\DeclareMathOperator{\im}{im}
\DeclareMathOperator{\res}{res}
\DeclareMathOperator{\ind}{ind}
\DeclareMathOperator{\algind}{algind}

\DeclareMathOperator{\id}{id}

\newenvironment{bnum}
{\begin{list}{}
    {\setlength{\labelwidth}{15pt}
     \setlength{\leftmargin}{\labelwidth}
    }
}
{\end{list}}

\numberwithin{equation}{section}
\theoremstyle{change}
\newtheorem{theorem}{Theorem}[section]
\newtheorem{prop}[theorem]{Proposition}
\newtheorem{lemma}[theorem]{Lemma}
\newtheorem{cor}[theorem]{Corollary}
\newtheorem{definition}[theorem]{Definition}

\begin{document}

\title[Free quantum groups]{The $ K $-theory of free quantum groups}

\author{Roland Vergnioux and Christian Voigt}

\address{Roland Vergnioux \\
         Laboratoire de Mathématiques Nicolas Oresme \\
         Université de Caen Basse-Normandie \\
         BP 5186 \\
         14032 Caen cedex \\
         France}
\email{roland.vergnioux@math.unicaen.fr}

\address{Christian Voigt \\
         Mathematisches Institut\\
         Westf\"alische Wilhelms-Universit\"at M\"unster \\
         Einsteinstra\ss e 62 \\
         48149 M\"unster\\
         Germany}
\email{cvoigt@math.uni-muenster.de}


\begin{abstract}
In this paper we study the $ K $-theory of free quantum groups in the sense of Wang and Van Daele, more precisely, of
free products of free unitary and free orthogonal quantum groups. We show that these quantum groups are $ K $-amenable 
and establish an analogue of the Pimsner-Voiculescu exact sequence. As a consequence, we obtain in particular an explicit 
computation of the $ K $-theory of free quantum groups. \\
Our approach relies on a generalization of methods from the Baum-Connes conjecture to the framework of discrete quantum 
groups. This is based on the categorical reformulation of the Baum-Connes conjecture developed by Meyer and Nest. As a 
main result we show that free quantum groups have a $ \gamma $-element and that $ \gamma = 1 $. \\
As an important ingredient in the proof we adapt the Dirac-dual Dirac method for groups acting on trees to the 
quantum case. We use this to extend some permanence properties of the Baum-Connes conjecture to our setting.
\end{abstract}

\maketitle

\section{Introduction} 

A classical result in the theory of $ C^* $-algebras is the computation of the $ K $-theory of the reduced 
group $ C^* $-algebra $ C^*_\red(\mathbb{F}_n) $ of the free group on $ n $ generators by Pimsner and Voiculescu \cite{PVfree}. 
Their result resolved in particular Kadison's problem on the existence of nontrivial projections in 
these $ C^* $-algebras. More generally, Pimsner and Voiculescu established an exact sequence for the $ K $-theory of 
reduced crossed products by free groups \cite{PVsequences}, \cite{PVfree}. This exact sequence is an important tool in 
operator $ K $-theory. \\
The $ K $-theory of the full group $ C^* $-algebra $ C^*_\max(\mathbb{F}_n) $ was calculated before by Cuntz in a simple and 
elegant way, based on a general formula for the $ K $-theory of free products \cite{Cuntzfreeproduct}. Motivated by this, 
Cuntz introduced the notion of $ K $-amenability for discrete groups and gave a shorter proof of the results of Pimsner and 
Voiculescu \cite{Cuntzkam}. The fact that free groups are $ K $-amenable expresses in a conceptually clear way that 
full and reduced crossed products for these groups cannot be distinguished on the level of $ K $-theory. \\
The main aim of this paper is to obtain analogous results for free quantum groups. In fact, in the theory of discrete quantum 
groups, the r\^ole of free quantum groups is analogous to the r\^ole of free groups among classical discrete groups. Roughly 
speaking, any discrete quantum group can be obtained as a quotient of a free quantum group. Classically, the free group on 
$ n $ generators can be described as the free product
$$
\mathbb{F}_n = \mathbb{Z} * \cdots * \mathbb{Z} 
$$
of $ n $ copies of $ \mathbb{Z} $. In the quantum case there is a similar free product construction, but in contrast to the 
classical situation there are different building blocks out of which free quantum groups are assembled. More precisely, a free 
quantum group is of the form 
$$
\mathbb{F}U(P_1) * \cdots * \mathbb{F}U(P_k) * \mathbb{F}O(Q_1) * \cdots * \mathbb{F}O(Q_l)
$$
for matrices $ P_i \in GL_{m_i}(\mathbb{C}) $ and $ Q_j \in GL_{n_j}(\mathbb{C}) $ such that $ \overline{Q_j} Q_j = \pm 1 $. Here 
we denote by $ \mathbb{F}U(P) $ and $ \mathbb{F}O(Q) $ the free unitary and free orthogonal quantum groups introduced by Wang 
and Van Daele \cite{Wangfree}, \cite{vDWuniversal}. The special case $ l = 0 $ and $ P_1 = \cdots = P_k = 1 \in GL_1(\mathbb{C}) $ 
of this family reduces to the classical free group $ \mathbb{F}_k $ on $ k $ generators. \\ 
In order to explain our notation let us briefly review some definitions. Given a matrix $ Q \in GL_n(\mathbb{C}) $, the full 
$ C^* $-algebra of the free unitary quantum group $ \mathbb{F}U(Q) $ is the universal $ C^* $-algebra $ C^*_\max(\mathbb{F}U(Q)) $ 
generated by the entries of an $ n \times n $-matrix $ u $ satisfying the relations that $ u $ and $ Q \overline{u} Q^{-1} $ are 
unitary. Here $ \overline{u} $ denotes the matrix obtained from $ u $ by taking the transpose of its adjoint. The full $ C^* $-algebra 
$ C^*_\max(\mathbb{F}O(Q)) $ of the free orthogonal quantum group $ \mathbb{F}O(Q) $ is the quotient of $ C^*_\max(\mathbb{F}U(Q)) $ 
by the relation $ u = Q \overline{u} Q^{-1} $. Finally, the full $ C^* $-algebra of the free product $ G * H $ of two discrete quantum 
groups $ G $ and $ H $ is the unital free product $ C^*_\max(G * H) = C^*_\max(G) * C^*_\max(H) $ of the corresponding full 
$ C^* $-algebras, see \cite{Wangfree}. That is, the full $ C^* $-algebra of a free quantum group as above is simply the free product 
of the $ C^* $-algebras $ C^*_\max(\mathbb{F}U(P_i)) $ and $ C^*_\max(\mathbb{F}O(Q_j)) $. \\
We remark that one usually writes $ C^*_\max(\mathbb{F}U(Q)) = A_u(Q) $ and $ C^*_\max(\mathbb{F}O(Q)) = A_o(Q) $ for 
these $ C^* $-algebras, compare \cite{Banicaunitary}. Following \cite{Voigtbcfo}, we use a different notation in order to emphasize 
that we shall view $ A_u(Q) $ and $ A_o(Q) $ as the full group $ C^* $-algebras of discrete quantum groups, and not as 
function algebras of compact quantum groups. \\
The approach to the $ K $-theory of free quantum groups in this paper is based on ideas and methods originating 
from the Baum-Connes conjecture \cite{BC}, \cite{BCH}. It relies in particular on the categorical 
reformulation of the Baum-Connes conjecture developed by Meyer and Nest \cite{MNtriangulated}. In fact, our 
main result is that free quantum groups satisfy an analogue of the strong Baum-Connes conjecture. The precise 
meaning of this statement will be explained in section \ref{secbc}, along with the necessary preparations from 
the theory of triangulated categories. Together with the results of Banica on the representation theory of 
free quantum groups \cite{Banicafo}, \cite{Banicaunitary}, the strong Baum-Connes property implies that every 
object of the equivariant Kasparov category $ KK^G $ for a free quantum group $ G $ has a projective resolution 
of length one. Based on this, we immediately obtain our Pimsner-Voiculescu exact sequence by invoking some 
general categorical considerations \cite{MNhomalg1}. As a consequence we conclude in particular that the reduced $ C^* $-algebras 
of unimodular free quantum groups do not contain nontrivial idempotents, extending the results of Pimsner and 
Voiculescu for free groups mentioned in the beginning. \\
This paper can be viewed as a continuation of \cite{Voigtbcfo}, where the Baum-Connes 
conjecture for free orthogonal quantum groups was studied. Our results here rely on the work in \cite{Voigtbcfo} 
on the one hand, and on geometric arguments using actions on quantum trees in the spirit of \cite{Vergniouxtrees} 
on the other hand. \\
To explain the general strategy let us consider first the case of free unitary quantum groups $ G = \mathbb{F}U(Q) $ 
for $ Q \in GL_n(\mathbb{C}) $ satisfying $ Q \overline{Q} = \pm 1 $. It was shown by Banica \cite{Banicaunitary} 
that $ \mathbb{F}U(Q) $ is a quantum subgroup of the free product $ \mathbb{F}O(Q) * \mathbb{Z} $ in this case. 
Since both $ \mathbb{F}O(Q) $ and $ \mathbb{Z} $ satisfy the strong Baum-Connes conjecture \cite{Voigtbcfo}, 
\cite{HKatmenable}, it suffices to prove inheritance properties of the conjecture for free products of quantum 
groups and for suitable quantum subgroups. In the case of free products we adapt the construction of Kasparov 
and Skandalis in \cite{KSbuildings} for groups acting on trees, and this is where certain quantum trees show up 
naturally. An important difference to the classical situation is that one has to work equivariantly with respect 
to the Drinfeld double $ \DD(G) $ of $ G $. \\
The quantum group $ \mathbb{F}U(Q) $ associated to a general matrix $ Q \in GL_n(\mathbb{C}) $ does not admit 
an embedding into a free product as above. As in \cite{Voigtbcfo} we use an indirect argument based on the
monoidal equivalences for free quantum groups obtained by Bichon, de Rijdt, and Vaes \cite{BdRV}. This allows 
us to reduce to matrices $ Q \in GL_2(\mathbb{C}) $, and in this case one may even assume $ Q \overline{Q} = \pm 1 $ 
without loss of generality. We might actually restrict attention to $ 2 \times 2 $-matrices throughout, however, this would 
not simplify the arguments. \\
Let us now describe how the paper is organized. In section \ref{secqg} we collect some preliminaries on quantum groups and fix 
our notation. Section \ref{secquantumsub} contains basic facts about quantum subgroups of discrete quantum groups 
and their homogeneous spaces. In section \ref{secdivisible} we introduce and discuss the notion of a divisible quantum subgroup
of a discrete quantum group. This concept appears naturally in the context of inheritance properties of the Baum-Connes conjecture. 
Roughly speaking, divisible quantum subgroups are particularly well-behaved from the point of view of corepresentation theory. In 
section \ref{secdirac} we define the Dirac element associated to a free product of discrete quantum groups acting on the corresponding 
quantum tree. Moreover we define the dual Dirac element and show that the resulting $ \gamma $-element is equal to the identity. 
In section \ref{secbc} we review the approach to the Baum-Connes conjecture developed by Meyer and Nest. We explain in particular 
the categorical ingredients needed to formulate the strong Baum-Connes property in our context. Then, using the considerations 
from section \ref{secdirac} and \cite{Voigtbcfo}, we prove that free quantum groups have the strong Baum-Connes property. Finally, in 
section \ref{secktheory} we discuss the main consequences of this result. As indicated above, we derive in particular the 
$ K $-amenability of free quantum groups and establish an analogue of the Pimsner-Voiculescu sequence. \\
Let us make some remarks on notation. We write $ \LH(\E) $ for the space of adjointable operators on a Hilbert $ A $-module $ \E $. 
Moreover $ \KH(\E) $ denotes the space of compact operators. The closed linear span of a subset $ X $ of a Banach space is denoted 
by $ [X] $. Depending on the context, the symbol $ \otimes $ denotes either the tensor product of Hilbert spaces, the minimal tensor 
product of $ C^* $-algebras, or the tensor product of von Neumann algebras. We write $ \odot $ for algebraic tensor products. For 
operators on multiple tensor products we use the leg numbering notation. If $ A $ and $ B $ are unital $ C^* $-algebras we 
write $ A * B $ for the unital free product, that is, the free product of $ A $ and $ B $ amalgamated over $ \mathbb{C} $. \\
We would like to thank G. Skandalis for helpful comments. 

\section{Preliminaries on quantum groups} \label{secqg} 

In this section we recall some basic definitions concerning quantum groups and their actions on $ C^* $-algebras. In 
addition we review the definition of the Drinfeld double and the description of its actions in terms of the underlying 
quantum group and its dual. For more information we refer to \cite{BSUM}, \cite{Kustermansuniversal}, 
\cite{KVLCQG}, \cite{Vaesimprimitivity}, \cite{Woronowiczleshouches}. Our notation and conventions will mainly follow 
\cite{NVpoincare}, \cite{Voigtbcfo}. \\
Let $ \phi $ be a normal, semifinite and faithful weight on a von Neumann algebra $ M $. We use the standard notation
$$
\M^+_\phi = \{ x \in M_+| \phi(x) < \infty \}, \qquad \N_\phi = \{ x \in M| \phi(x^* x) < \infty \}
$$
and write $ M_*^+ $ for the space of positive normal linear functionals on $ M $. If 
$ \Delta: M \rightarrow M \otimes M $ is a normal unital $ * $-homomorphism, the weight $ \phi $ is called left invariant
with respect to $ \Delta $ provided
$$
\phi((\omega \otimes \id)\Delta(x)) = \phi(x) \omega(1)
$$
for all $ x \in \M_\phi^+ $ and $ \omega \in M_*^+ $. Similarly one defines right invariant weights.
\begin{definition} \label{defqg}
A locally compact quantum group $ G $ is given by a von Neumann algebra $ L^\infty(G) $ together with a normal unital $ * $-homomorphism
$ \Delta: L^\infty(G) \rightarrow L^\infty(G) \otimes L^\infty(G) $ satisfying the coassociativity relation
$$
(\Delta \otimes \id)\Delta = (\id \otimes \Delta)\Delta
$$
and normal semifinite faithful weights $ \phi $ and $ \psi $ on $ L^\infty(G) $ which are left and right invariant, respectively.
\end{definition}
If $ G $ is a locally compact group, then the algebra $ L^\infty(G) $ of essentially bounded measurable functions on $ G $ together with the comultiplication
$ \Delta: L^\infty(G) \rightarrow L^\infty(G) \otimes L^\infty(G) $ given by
$$
\Delta(f)(s,t) = f(st)
$$
defines a locally compact quantum group. The weights $ \phi $ and $ \psi $ are given in this case by left and right Haar measures 
on $ G $, respectively. Of course, for a general locally compact quantum group $ G $ the notation $ L^\infty(G) $ is purely formal. \\
An important tool in the study of locally compact quantum groups are multiplicative unitaries. If $ \Lambda: \N_\phi \rightarrow \HH_G = L^2(G) $ 
is a GNS-construction for the weight $ \phi $, then the multiplicative unitary $ W_G = W $ is the operator on $ \HH_G \otimes \HH_G $ given by 
$$
W^*(\Lambda(f) \otimes \Lambda(g)) = (\Lambda \otimes \Lambda)(\Delta(g)(f \otimes 1))
$$
for all $ f, g \in \N_\phi $. This unitary satisfies the pentagonal equation
$
W_{12} W_{13} W_{23} = W_{23} W_{12}, 
$
and one can recover the von Neumann algebra $ L^\infty(G) $ as the strong closure of the algebra 
$ (\id \otimes \LH(\HH_G)_*)(W) $ where $ \LH(\HH_G)_* $ denotes the space of normal linear functionals on $ \LH(\HH_G) $. 
Moreover one has
$$
\Delta(f) = W^*(1 \otimes f) W
$$
for all $ f \in L^\infty(G) $. Let us remark that we will only consider quantum groups $ G $ for which $ \HH_G = L^2(G) $ is a separable 
Hilbert space. \\
The group-von Neumann algebra $ \mathcal{L}(G) $ of the quantum group $ G $ is the strong
closure of the algebra $ (\LH(\HH_G)_* \otimes \id)(W) $ with the comultiplication 
$ \hat{\Delta}: \mathcal{L}(G) \rightarrow \mathcal{L}(G) \otimes \mathcal{L}(G) $ given by
$$
\hat{\Delta}(x) = \hat{W}^* (1 \otimes x) \hat{W}
$$
where $ \hat{W} = \Sigma W^* \Sigma $ and $ \Sigma \in \LH(\HH_G \otimes \HH_G) $ is the flip map. It defines
a locally compact quantum group $ \hat{G} $ which is called the dual of $ G $. The left invariant weight
$ \hat{\phi} $ for the dual quantum group has a GNS-construction $ \hat{\Lambda}: \N_{\hat{\phi}} \rightarrow \HH_G $,
and according to our conventions we have $ \mathcal{L}(G) = L^\infty(\hat{G}) $. \\
The reduced $ C^* $-algebra of functions on the quantum group $ G $ is 
$$
C^\red_0(G) = [(\id \otimes \LH(\HH_G)_*)(W)] \subset L^\infty(G), 
$$
and the reduced group $ C^* $-algebra of $ G $ is 
$$
C^*_\red(G) = [(\LH(\HH_G)_* \otimes \id)(W)] \subset \L(G).  
$$
Moreover we have $ W \in M(C_0^\red(G) \otimes C^*_\red(G)) $. 
Together with the comultiplications inherited from $ L^\infty(G) $ and $ \mathcal{L}(G) $, respectively,
the $ C^* $-algebras $ C^\red_0(G) $ and $ C^*_\red(G) $ are Hopf-$ C^* $-algebras in the following sense.
\begin{definition}\label{defhopfcstar}
A Hopf $ C^* $-algebra is a $ C^* $-algebra $ S $ together with an injective nondegenerate $ * $-homomorphism
$ \Delta: S \rightarrow M(S \otimes S) $ such that the diagram
$$
\xymatrix{
S \ar@{->}[r]^{\Delta} \ar@{->}[d] & M(S \otimes S) \ar@{->}[d]^{\id \otimes \Delta} \\
M(S \otimes S) \ar@{->}[r]^{\!\!\!\!\!\!\!\!\! \Delta \otimes \id} & M(S \otimes S \otimes S)
     }
$$
is commutative and $ [\Delta(S)(1 \otimes S)] = S \otimes S = [(S \otimes 1)\Delta(S)] $. \\
A morphism between Hopf-$ C^* $-algebras $ (S, \Delta_S) $ and $ (T, \Delta_T) $ is a nondegenerate $ * $-homomorphism
$ \pi: S \rightarrow M(T) $ such that $ \Delta_T\, \pi = (\pi \otimes \pi)\Delta_S $.
\end{definition}
If $ S $ is a Hopf $ C^* $-algebra we write $ S^\cop $ for $ S $ viewed as a Hopf $ C^* $-algebra with the opposite coproduct 
$ \Delta^\cop = \sigma \Delta $, where $ \sigma $ denotes the flip map. \\
We are mainly interested in discrete or compact quantum groups, and in these cases the general theory simplifies considerably. 
A quantum $ G $ is called compact if $ C^\red_0(G) $ is unital, and it is called discrete if $ C^*_\red(G) $ is unital. 
We note that $ G $ is compact iff $ \hat{G} $ is discrete, and vice versa. For a discrete quantum group $ G $ we 
will write $ l^2(G) = \HH_G $ for the Hilbert space associated to $ G $. 
\begin{definition} 
Let $ B $ be a $ C^* $-algebra. A unitary corepresentation of a Hopf-$ C^* $-algebra $ S $ on a Hilbert $ B $-module $ \E $ 
is a unitary $ X \in \LH(S \otimes \E) $ satisfying
$$
(\Delta \otimes \id)(X) = X_{13} X_{23}.
$$
\end{definition} 
If $ S = C^*_\red(G) $ for a locally compact quantum group $ G $ and $ X $ is a unitary corepresentation of $ S $ on 
$ \E $ we say that $ \E $ is a unitary corepresentation of $ G $. If $ G $ is discrete, all unitary corepresentations
of $ G $ on Hilbert spaces are completely reducible, and all irreducible corepresentations are finite dimensional. 
We will write $ \Corep(G) $ for the corresponding semisimple $ C^* $-tensor category of finite dimensional corepresentations. 
Moreover we denote by $ \Irr(G) $ the set of equivalence classes of irreducible corepresentations of $ G $ in
this case. \\
If $ S $ is a Hopf-$ C^* $-algebra, then a universal dual of $ S $ is a Hopf-$ C^* $-algebra $ \hat{S} $ together with a unitary
corepresentation $ \X \in M(S \otimes \hat{S}) $ satisfying the following universal property.
For every Hilbert $ B $-module $ \E $ and every unitary corepresentation $ X \in \LH(S \otimes \E) $ there exists a unique 
nondegenerate $ * $-homomorphism $ \pi_X: \hat{S} \rightarrow \LH(\E) $ such that $ (\id \otimes \pi_X)(\X) = X $. \\
Every locally compact quantum group $ G $ admits a universal dual $ C^*_\max(G) $ of $ C_0^\red(G) $ and
a universal dual $ C^\max_0(G) $ of $ C^*_\red(G) $, respectively \cite{Kustermansuniversal}.
In general, we have a surjective morphism $ \hat{\pi}: C^*_\max(G) \rightarrow C^*_\red(G) $ of Hopf-$ C^* $-algebras
associated to the multiplicative unitary $ W \in M(C_0(G) \otimes C^*_\red(G)) $. Similarly, there is 
a canonical surjective morphism $ \pi: C^\max_0(G) \rightarrow C^\red_0(G) $. Every discrete 
quantum group $ G $ is coamenable in the sense that $ \pi $ is an isomorphism. We will 
simply write $ C_0(G) $ for $ C^\max_0(G) = C^\red_0(G) $ in this case. \\
If $ G $ is a discrete quantum group, then inside $ C^*_\red(G) $ we have the dense Hopf-$ * $-algebra $ \mathbb{C}[G] $ 
spanned by the matrix coefficients of all finite dimensional unitary corepresentations of $ C^*_\red(G) $, and 
inside $ C_0(G) $ we have a dense multiplier Hopf-$ * $-algebra $ C_c(G) $ in the sense of Van Daele \cite{vDadvances}. 
More precisely, the space of matrix coefficients $ \mathbb{C}[G] $ can be identified with a dense linear subspace $ E $ of $ l^2(G) $, 
and the algebras $ \mathbb{C}[G] \subset C^*_\red(G) $ and $ C_c(G) \subset C_0(G) $ are obtained from the elements 
$ \omega \in \LH(\HH_G)_* $ associated to vectors from $ E \subset \HH_G = l^2(G) $. 
If $ f \in C_c(G) $ and $ x \in \mathbb{C}[G] $ are represented by $ L_f, L_x \in \LH(l^2(G))_* $ in the sense 
that $ (\id \otimes L_f)(W) = f $ and $ (L_x \otimes \id)(W) = x $, 
then we obtain a well-defined bilinear pairing 
\begin{equation*}
\bra f, x \ket = \bra x, f \ket = (L_x \otimes L_f)(W) = L_f(x) = L_x(f)
\end{equation*}
between $ C_c(G) $ and $ \CH[G] $, see definition 1.3 of \cite{BSUM}. 
It is easy to check that the product of $ \CH[G] $ is dual to the coproduct
of $ C_c(G) $, whereas the product of $ C_c(G) $ is dual to the opposite 
coproduct of $ \CH[G] $. In other terms, we have for all $ f, g \in C_c(G)$ and $ x, y \in \CH[G] $ 
the relations 
\begin{equation*}
\bra f, xy \ket = \bra f_{(1)}, x \ket \bra f_{(2)}, y \ket \quad \text{and} \quad
\bra fg, x \ket = \bra f, x_{(2)} \ket \bra g, x_{(1)} \ket
\end{equation*}
where we use the Sweedler notation $ \Delta(f) = f_{(1)} \odot f_{(2)} $ and $ \hat{\Delta}(x) = x_{(1)} \odot x_{(2)}$
for the comultiplications on $ C_c(G)$ and $ \CH[G] $ induced by $ W $ as above. We point out that 
this notation has to be interpreted with care, let us remark in particular that the coproduct $ \Delta(f) $ of an element 
$ f $ of the multiplier Hopf $ * $-algebra $ C_c(G) $ can be represented only as an infinite sum of simple tensors in general. \\
At several points of the paper we will consider free products of discrete quantum groups. 
If $ G $ and $ G $ are discrete quantum groups the free product $ G * H $ is the discrete 
quantum group determined by $ C^*_\max(G * H) = C^*_\max(G) * C^*_\max(H) $, equipped with the comultiplication 
induced from the two factors in the evident way, see \cite{Wangfree} for more information. 
The irreducible corepresentations of $ G * H $ are precisely the alternating tensor
products of nontrivial irreducible corepresentations of $ G $ and $ H $. We may 
therefore identify $ \Irr(G * H) $ with the set $ \Irr(G) * \Irr(H) $ of alternating words 
in $ \Irr(G) $ and $ \Irr(H) $. \\
We are mainly interested in the free unitary and free orthogonal quantum groups introduced by Wang and 
Van Daele \cite{Wangfree}, \cite{vDWuniversal}. These discrete quantum groups are most conveniently 
defined in terms of their full group $ C^* $-algebras. 
For a matrix $ u = (u_{ij}) $ of elements in a $ * $-algebra we write
$ \overline{u} = (u_{ij}^*) $ and $ u^t = (u_{ji}) $ for its conjugate and transposed matrices, respectively. 
\begin{definition} \label{deffu}
Let $ n \in \mathbb{N} $ and $ Q \in GL_n(\mathbb{C}) $. The group $ C^* $-algebra $ C^*_\max(\mathbb{F}U(Q)) $ 
of the free unitary quantum group $ \mathbb{F}U(Q) $ is the universal $ C^* $-algebra with generators 
$ u_{ij}, 1 \leq i,j \leq n $ such that the resulting matrices $ u $ and $ Q \overline{u} Q^{-1} $ are unitary. 
The comultiplication $ \hat{\Delta}: C^*_\max(\mathbb{F}U(Q)) \rightarrow C^*_\max(\mathbb{F}U(Q)) \otimes C^*_\max(\mathbb{F}U(Q)) $ 
is given by 
$$
\hat{\Delta}(u_{ij}) = \sum_{k = 1}^n u_{ik} \otimes u_{kj}
$$
on the generators. 
\end{definition} 
As explained in the introduction, we adopt the conventions in \cite{Voigtbcfo} and deviate from the standard notation 
$ A_u(Q) $ for the $ C^* $-algebras $ C^*_\max(\mathbb{F}U(Q)) $. Let us also remark that there is a canonical isomorphism 
$ C^*_\max(\mathbb{F}U(Q))^\cop \cong C^*_\max(\mathbb{F}U(Q^t)) $ for all $ Q \in GL_n(\mathbb{C}) $. 
\begin{definition} \label{deffo}
Let $ Q \in GL_n(\mathbb{C}) $ such that $ Q \overline{Q} = \pm 1 $. The group $ C^* $-algebra $ C^*_\max(\mathbb{F}O(Q)) $ 
of the free orthogonal quantum group $ \mathbb{F}O(Q) $ is the universal $ C^* $-algebra with generators 
$ u_{ij}, 1 \leq i,j \leq n $ such that the resulting matrix $ u $ is unitary and the relation$ u = Q \overline{u} Q^{-1} $ 
holds. The comultiplication of the generators is given by the same formula as in the unitary case.
\end{definition} 
We remark that the free quantum groups $ \mathbb{F}O(Q) $ for $ Q \in GL_2(\mathbb{C}) $ exhaust up to isomorphism precisely 
the duals of $ SU_q(2) $ for $ q \in [-1,1] \setminus \{0\} $. \\
If $ Q \in GL_n(\mathbb{C}) $ is arbitrary the above definition of $ \mathbb{F}O(Q) $ still makes sense. However, in this 
case the fundamental corepresentation $ u $ is no longer irreducible in general. According to \cite{Wangstructure}, the quantum group 
$ \mathbb{F}O(Q) $ can be decomposed into a free product of the form 
$$
\mathbb{F}O(Q) \cong \mathbb{F}U(P_1) * \cdots * \mathbb{F}U(P_k) * \mathbb{F}O(Q_1) * \cdots * \mathbb{F}O(Q_l)
$$ 
for appropriate matrices $ P_i $ and $ Q_j $ such that $ Q_j \overline{Q_j} = \pm 1 $ for all $ j $. In this way the study of 
$ \mathbb{F}O(Q) $ for general $ Q $ reduces to the above cases. \\
Let us now fix our notation concerning coactions on $ C^* $-algebras and crossed products. 
\begin{definition} \label{defcoaction}
A (continuous, left) coaction of a Hopf $ C^* $-algebra $ S $ on a $ C^* $-algebra $ A $ is an injective nondegenerate $ * $-homomorphism
$ \alpha: A \rightarrow M(S \otimes A) $ such that the diagram
$$
\xymatrix{
A \ar@{->}[r]^{\alpha} \ar@{->}[d]^\alpha & M(S \otimes A) \ar@{->}[d]^{\Delta \otimes \id} \\
M(S \otimes A) \ar@{->}[r]^{\!\!\!\!\!\!\!\!\!\! \id \otimes \alpha} & M(S \otimes S \otimes A)
     }
$$
is commutative and $ [\alpha(A)(S \otimes 1)] = S \otimes A $. \\
If $ (A, \alpha) $ and $ (B, \beta) $ are $ C^* $-algebras with coactions of $ S $, then a $ * $-homomorphism 
$ f: A \rightarrow M(B) $ is called $ S $-colinear, or $ S $-equivariant, if $ \beta f = (\id \otimes f)\alpha $.
\end{definition}
A $ C^* $-algebra $ A $ equipped with a continuous coaction of the Hopf-$ C^* $-algebra $ S $ is called an $ S $-$ C^* $-algebra.
If $ S = C^\red_0(G) $ for a locally compact quantum group $ G $ we say that $ A $ is $ G $-$ C^* $-algebra. In this case 
$ S $-colinear $ * $-homomorphisms will be called $ G $-equivariant or simply equivariant. \\
If $ G $ is a discrete quantum group it is useful to consider algebraic coactions as well. 
Assume that $ A $ is a $ * $-algebra equipped with an injective $ * $-homomorphism $ \alpha: A \rightarrow \M(C_c(G) \odot A) $ such 
that $ (C_c(G) \odot 1)\alpha(A) = C_c(G) \odot A $ and $ (\Delta \odot \id)\alpha = (\id \odot \alpha)\alpha $. 
Here $ \M(C_c(G) \odot A) $ is the algebraic multiplier algebra of $ C_c(G) \odot A $, see \cite{vDadvances}. 
In this case we say that $ A $ is a $ G $-algebra, and we refer to $ \alpha $ as an algebraic coaction. \\
If $ A $ is a $ G $-algebra we have a left $ \CH[G] $-module structure on $ A $ given by 
$$ 
x \cdot a = \bra a_{(-1)}, S(x) \ket a_{(0)} 
$$ 
where $ \alpha(a) = a_{(-1)} \odot a_{(0)} $ is Sweedler notation for the coaction on $ A $ and $ S $ is the antipode of $ \CH[G] $. 
This action turns $ A $ into a $ \CH[G] $-module algebra in the sense that 
$$ 
x \cdot (ab) = (x_{(1)} \cdot a)(x_{(2)} \cdot b), \qquad x \cdot a^* = (S(x)^* \cdot a)^*  
$$ 
for all $ x \in \CH[G] $ and $ a, b \in A $. \\
Similarly, we may associate to an algebraic coaction $ \gamma: A \rightarrow \M(\CH[G] \odot A) $ a left 
$ C_c(G) $-module structure on $ A $ given by 
$$ 
f \cdot a = \bra f, a_{(-1)} \ket a_{(0)} 
$$ 
where $ \gamma(a) = a_{(-1)} \odot a_{(0)} $. This turns $ A $ into a $ C_c(G) $-module algebra. We will study 
coactions in this way at several points below. \\
Let us next recall the definition of reduced crossed products. If $ G $ is a locally compact quantum group and $ A $ is a 
$ G $-$ C^* $-algebra with coaction $ \alpha $, the reduced crossed product $ G \ltimes_\red A = C^*_\red(G)^\cop \ltimes_\red A $ is defined by 
\begin{displaymath}
G \ltimes_\red A = C^*_\red(G)^\cop \ltimes_\red A = [(C^*_\red(G) \otimes 1) \alpha(A)]  \subset M(\KH(l^2(G)) \otimes A).
\end{displaymath}
The crossed product $ G \ltimes_\red A $ is naturally a $ C^*_\red(G)^\cop $-$ C^* $-algebra with the dual coaction 
$ \hat{\alpha}: G \ltimes_\red A \rightarrow M(C^*_\red(G)^\cop \otimes (G \ltimes_\red A))$ determined by 
$ \hat{\alpha}(x\otimes 1) = \hat\Delta^\cop(x) \otimes 1 $ for $ x \in C^*_\red(G) $ and 
$ \hat{\alpha}(\alpha(a)) = 1 \otimes \alpha(a) $ for $ a \in A $. \\
Apart from the reduced crossed product $ G \ltimes_\red A $ one also has the full crossed product 
$ G \ltimes_\max A = C^*_\max(G)^\cop \ltimes_\max A $ which is defined by a universal property. There is a canonical map 
$ G \ltimes_\max A \rightarrow G \ltimes_\red A $, and this map is an isomorphism if $ G $ is amenable. We refer 
to \cite{NVpoincare} for more details. \\  
Finally, let us recall the definition of the Drinfeld double of a locally compact quantum group $ G $, see \cite{BV}. 
The Drinfeld double $ \DD(G) $ of $ G $ is a locally compact quantum group such that the Hopf $ C^* $-algebra 
$ C^\red_0(\DD(G_q)) $ is given by $ C^\red_0(\DD(G_q)) = C_0^\red(G) \otimes C^*_\red(G) $ with the comultiplication
$$
\Delta_{\DD(G_q)} = (\id \otimes \sigma \otimes \id)(\id \otimes \ad(W) \otimes \id)(\Delta \otimes \hat{\Delta}), 
$$
here $ \ad(W) $ is conjugation with the multiplicative unitary $ W \in M(C_0^\red(G) \otimes C^*_\red(G)) $ 
and $ \sigma $ denotes the flip map. \\
It is shown in \cite{NVpoincare} that a $ \DD(G) $-$ C^* $-algebra $ A $ is uniquely determined by coactions 
$ \alpha: A \rightarrow M(C_0^\red(G) \otimes A) $ and $ \gamma: A \rightarrow M(C^*_\red(G) \otimes A) $ satisfying the 
Yetter-Drinfeld compatibility condition
$$
(\sigma \otimes \id)(\id \otimes \alpha)\gamma = (\ad(W) \otimes \id) (\id \otimes \gamma)\alpha.
$$
In a similar way on can study $ \DD(G) $-equivariant Hilbert modules. \\
The Drinfeld double plays an important r\^ole in the definition of braided tensor products. 
If $ A $ is a $ \DD(G) $-$ C^* $-algebra determined by the coactions $ \alpha: A \rightarrow M(C_0^\red(G) \otimes A) $ and 
$ \gamma: A \rightarrow M(C^*_\red(G) \otimes A) $ and $ B $ is $ G $-$ C^* $-algebra with coaction 
$ \beta: B \rightarrow M(C_0^\red(G) \otimes B) $, then the braided tensor product of $ A $ and $ B $ is defined by 
$$
A \twisted B = A \twisted_G B = [\gamma(A)_{12} \beta(B)_{13}] \subset \LH(\HH_G \otimes A \otimes B). 
$$
The braided tensor product is naturally a $ G $-$ C^* $-algebra, and it is a natural replacement of the minimal tensor product 
of $ G $-$ C^* $-algebras in the group case, see \cite{NVpoincare} for more information. \\
We will be interested in the case that $ G $ is a discrete quantum group and work with $ \DD(G) $-$ C^* $-algebras obtained from 
algebraic actions and coactions. More precisely, in the algebraic setting we have a pair of algebraic coactions 
$ \alpha: A \rightarrow \M(C_c(G) \odot A) $ and $ \gamma: A \rightarrow \M(\CH[G] \odot A) $, and the Yetter-Drinfeld condition can be 
written as
$$
f_{(1)} a_{(-1)} S(f_{(3)}) \odot f_{(2)} \cdot a_{(0)} = (f \cdot a)_{(-1)} \odot (f \cdot a)_{(0)}
$$
in this case. Here we denote by $ f \cdot a $ the action of $ f \in C_c(G) $ on $ a \in A $ corresponding to the coaction $ \gamma $ and 
$ \alpha(a) = a_{(-1)} \odot a_{(0)} $.

\section{Quantum subgroups of discrete quantum groups} \label{secquantumsub}

In this section we collect some basic constructions related to quantum subgroups of discrete quantum groups. The general 
concept of a closed quantum subgroup of a locally compact quantum group is discussed in \cite{Vaesimprimitivity}, 
\cite{VaesVainermanlowdim}. \\
A morphism $ H \rightarrow G $ of locally compact quantum groups is a nondegenerate $ * $-homomorphism 
$ \pi: C_0^\max(G) \rightarrow M(C_0^\max(H)) $ compatible with the comultiplications. For every such morphism there exists 
a dual morphism $ \hat{\pi}: C^*_\max(H) \rightarrow M(C^*_\max(G)) $. By definition, a closed quantum subgroup $ H \subset G $ 
of a locally compact quantum group $ G $ is a morphism $ H \rightarrow G $ for which the associated map 
$ C^*_\max(H) \rightarrow M(C^*_\max(G)) $ is accompanied by a compatible unital faithful normal $ * $-homomorphism 
$ \L(H) \rightarrow \L(G) $ of the von Neumann algebras \cite{Vaesimprimitivity}, \cite{VaesVainermanlowdim}. \\
In the case of discrete quantum groups we will simply speak of quantum subgroups instead of closed quantum subgroups. 
In fact, in this case the above definition can be rephrased purely algebraically. Let $ G $ be a discrete quantum group 
and let $ H \subset G $ be closed quantum subgroup. Then $ H $ is automatically discrete, and the corresponding map 
$ \iota: \L(H) \rightarrow \L(G) $ induces an injective $ * $-homomorphism $ \mathbb{C}[H] \rightarrow \mathbb{C}[G] $ of 
Hopf-$ * $-algebras. Conversely, if $ H $ and $ G $ are discrete, every injective homomorphism 
$ \iota: \mathbb{C}[H] \rightarrow \mathbb{C}[G] $ of Hopf-$ * $-algebras extends uniquely to a unital faithful normal 
$ * $-homomorphism on the von Neumann algebra level, and thus turns $ H $ into a quantum subgroup of $ G $. \\
If $ H \subset G $ is a quantum subgroup of the discrete quantum group $ G $, then the category $ \Corep(H) $ 
is a full tensor subcategory of $ \Corep(G) $ containing the trivial corepresentation and closed under taking duals. 
In fact, quantum subgroups of $ G $ can be characterized in terms of such subcategories, see for instance \cite{Vergniouxkam}. 
The $ C^* $-algebra $ C_0(H) $ is obtained as the sum of matrix blocks in $ C_0(G) $ corresponding to 
corepresentations in $ \Irr(H) $, and the map $ \pi: C_0(G) \rightarrow C_0(H) \subset M(C_0(H)) $
is the canonical projection. \\
We will need some basic facts concerning homogeneous spaces and induced $ C^* $-algebras for discrete quantum groups. 
The general theory of induced coactions for locally compact quantum groups is due to Vaes \cite{Vaesimprimitivity}, 
and it is technically quite involved. However, in the discrete case the constructions that we need can be described more 
directly. We start with some elementary algebraic considerations. \\
For a discrete quantum group $ G $ we denote by $ C(G) $ the algebra $ \M(C_c(G)) $ of algebraic multipliers of $ C_c(G) $. 
Heuristically, the elements of $ C(G) $ can be viewed as functions with arbitrary support on $ G $. 
The canonical bilinear pairing between $ C_c(G) $ and $ \mathbb{C}[G] $ extends uniquely to a bilinear pairing between 
$ C(G) $ and $ C_c(G) $, and in this way we obtain a natural identification of $ C(G) $ with the algebraic dual space 
$ \CH[G]^* $. We write again $ \bra f, x \ket = \bra x, f \ket $ for $ f \in C_c(G) $ and $ x \in \mathbb{C}[G] $. 
\begin{definition}
Let $ H $ be a quantum subgroup of a discrete quantum group $ G $. The algebra of all functions on the homogeneous 
space $ G/H $ is 
\begin{equation*}
C(G/H) = \{f \in C(G) \mid (\id \odot \pi)\Delta(f) = f \odot 1 \}
\end{equation*}
where $ \pi: C(G) \rightarrow C(H) $ is the canonical projection map. 
Moreover we define 
\begin{equation*}
\CH[G/H] = \CH[G] \odot_{\CH[H]}\CH 
\end{equation*}
where $ \CH[H] $ acts via $ \iota: \CH[H] \rightarrow \CH[G] $ on $ \CH[G] $ and via the counit on $ \CH $. 
\end{definition}
Let $ \tau: \CH[G] \rightarrow \CH[G/H] $ be the canonical projection given by $ \tau(x) = x \odot 1 $. The transpose of 
$ \tau $ induces an injective linear map $ \mathbb{C}[G/H]^* \rightarrow \mathbb{C}[G]^* = C(G) $. Moreover $ \pi $ is the 
transpose of $ \iota $ in the duality between $ C(G) $ and $ \CH[G] $. Hence for $ f \in C(G/H), x \in \CH[G] $ and 
$ y \in \CH[H] $ we obtain
\begin{align*}
\bra f, x \iota(y) \ket &= \bra \Delta(f), x \odot \iota(y) \ket = \bra (\id \odot \pi)\Delta(f), x \odot y \ket 
= \bra f \odot 1, x \odot y \ket = \bra f, x \epsilon(y) \ket
\end{align*}
and conclude that $ f $ is contained in $ \mathbb{C}[G/H]^* $. An analogous argument shows that every $ f \in \mathbb{C}[G/H]^* $ 
satisfies the invariance condition defining $ C(G/H) $. In other words, we may identify $ \mathbb{C}[G/H]^* \subset \mathbb{C}[G]^* $ 
with $ C(G/H) \subset C(G) $. Note also that the comultiplication of $ \CH[G] $ induces a natural coalgebra structure on $ \CH[G/H] $. \\
Let $ \Corep(H) \subset \Corep(G) $ be the full tensor subcategory corresponding to the inclusion $ H \subset G $. For an irreducible 
corepresentation $ r \in \Irr(G) $ we denote by $ \mathbb{C}[G]_r \subset \mathbb{C}[G] $ the linear subspace of matrix 
coefficients of $ r $. By definition of the tensor product we have 
$$ 
\mathbb{C}[G]_r \cdot \mathbb{C}[G]_s = \bigoplus_{t \subset r \otimes s} \mathbb{C}[G]_t 
$$ 
for $ r, s \in \Irr(G) $. Moreover if $ r \otimes s = t $ is already irreducible, then the multiplication map defines an isomorphism 
$ \mathbb{C}[G]_r \odot \mathbb{C}[G]_s \cong \mathbb{C}[G]_t $. \\
In \cite{Vergniouxkam}, an equivalence relation on $ \Irr(G) $ is defined by setting $ r \sim s $ iff $ \overline{r} \otimes s $ 
contains an element $ t \in \Irr(H) $. Equivalently, $ r \sim s $ iff $ s \subset r \otimes t $ for some $ t \in \Irr(H) $. 
We denote by $ \Irr(G)/\Irr(H) $ the corresponding quotient space, so that the class $ [r] $ of $ r \in \Irr(G) $ is the set of 
isomorphism classes of irreducible subobjects of tensor products $ r \otimes t $ with $ t \in \Irr(H) $. \\
For $ \alpha \in \Irr(G)/\Irr(H) $ let us write $ \mathbb{C}[G]_\alpha $ the direct sum of subspaces $ \mathbb{C}[G]_r $ with 
$ r \in \alpha $. Then $ \mathbb{C}[G] $ is clearly the direct sum of the subspaces $ \mathbb{C}[G]_\alpha $ over all 
$ \alpha \in \Irr(G)/\Irr(H) $, and $ \mathbb{C}[G]_{[\epsilon]} = \mathbb{C}[H] $ where $ \epsilon $ is the trivial 
corepresentation. \\
We may write $ C(G) $ as 
$$
C(G) \cong \prod_{r \in \Irr(G)} C(G)_r
$$ 
where $ C(G)_r = \mathbb{C}[G]_r^* $ is the matrix algebra corresponding to $ r \in \Irr(G) $. Let us denote by $ p_r \in C(G) $ the 
central projection associated to the identity in $ C(G)_r $, so that $ \bra p_r, x \ket = \delta_{t,r} \epsilon(x) $ for $ x \in \CH[G]_t $. 
For $ \alpha \in \Irr(G)/\Irr(H)$ we let $ p_\alpha \in C(G) $ be the sum of the projections $ p_r $ with $ r \in \alpha $. 
By construction, the sum of the elements $ p_\alpha $ over $ \alpha \in \Irr(G)/\Irr(H) $ is equal to $ 1 $. 
For $ \alpha, \beta \in \Irr(G)/\Irr(H)$ and $ x \in \CH[G]_\beta $ we compute 
$ \bra p_\alpha, x \ket = \delta_{\alpha, \beta} \epsilon(x) $ and conclude that 
$$
\bra p_\alpha, xy \ket = \bra p_\alpha, x \epsilon(y) \ket 
$$ 
for all $ x \in \CH[G] $ and $ y \in \CH[H] $. This implies $ p_\alpha \in C(G/H) $. \\
Let $ A $ be an $ H $-algebra with coaction $ \alpha: A \rightarrow \M(C_c(H) \odot A) $ in the sense explained in section \ref{secqg}. 
Inside the algebraic multiplier algebra $ \M(C_c(G) \odot A) $ we have the $ * $-subalgebra 
$$
C(G, A) = \prod_{r \in \Irr(G)} C(G)_r \odot A \subset \M(C_c(G) \odot A)
$$
of functions with values in $ A \subset \M(A) $. We let 
$$
C(G, A)^H = \{f \in C(G,A) \mid (\rho \odot \id)(f) = (\id \odot \alpha)(f) \}
$$
where $ \rho: C_c(G) \rightarrow \M(C_c(G) \odot C_c(H)) $ is the right coaction $ \rho = (\id \odot \pi)\Delta $. 
Note that $ p_\alpha \odot 1 $ is a central element in $ \M(C(G, A)^H) $ for every $ \alpha \in \Irr(G)/\Irr(H) $ in a natural way. 
In the case $ A = \mathbb{C} $ the construction of $ C(G,A)^H $ obviously reduces to the algebraic homogeneous space $ C(G/H) $ defined above. 
\begin{definition} 
Let $ G $ be a discrete quantum group and let $ H \subset G $ be a quantum subgroup. The algebraic induced algebra of an $ H $-algebra $ A $ 
is the $ * $-algebra
$$
\algind_H^G(A) = \bigoplus_{\alpha \in \Irr(G)/\Irr(H)} p_\alpha C(G, A)^H 
$$
inside $ C(G,A)^H $. 
\end{definition}
In the case $ A = \mathbb{C} $ with the trivial action we shall use the notation 
$$ 
C_c(G/H) = \bigoplus_{\alpha \in \Irr(G)/\Irr(H)} p_\alpha C(G/H) \subset C(G) 
$$
for the algebraic induced algebra. We may view $ C_c(G/H) $ as the algebra of finitely supported functions on the 
homogeneous space $ G/H $. \\
Let us show that $ \algind_H^G(A) $ becomes a $ G $-algebra in a natural way. To this end let $ r \in \Irr(G), \alpha \in \Irr(G)/\Irr(H) $ 
and consider the finite subset $ \overline{r} \alpha \subset \Irr(G)/\Irr(H) $ given by the equivalence classes of all irreducible subobjects 
of $ \overline{r} \otimes s $ for some $ s \in \alpha $. Note that this set is independent of the choice of $ s $. Moreover let 
$ p_{\overline{r} \alpha} \in C_c(G/H) $ be the sum of the projections $ p_\beta $ for $ \beta \in \overline{r} \alpha $. 
Then we have 
$$
(p_r \odot \id) \Delta(p_\alpha) = (p_r \odot p_{\overline{r} \alpha}) \Delta(p_\alpha), 
$$
and it follows that $ (C_c(G) \odot 1) (\Delta \odot \id)(\algind_H^G(A)) \subset C_c(G) \odot \algind_H^G(A) $ 
inside $ \M(C_c(G) \odot C_c(G) \odot A) $. Using the antipode of $ C_c(G) $ it is straightforard to check that this inclusion is 
in fact an equality. We conclude that the map 
$$ 
\Delta \odot \id: \algind_H^G(A) \rightarrow 
\M(C_c(G) \odot \algind_H^G(A)) 
$$ 
is a well-defined algebraic coaction which turns $ \algind_H^G(A) $ into a $ G $-algebra. \\
We are mainly interested in the case that $ A $ is an $ H $-$ C^* $-algebra. In this situation we let 
$ C_b(G, A) \subset C(G, A) $ be the $ l^\infty $-direct sum of the spaces $ C(G)_t \otimes A $ over $ t \in \Irr(G) $. Note 
that $ C_b(G, A) $ is naturally contained in the multiplier algebra $ M(C_0(G) \otimes A) $. 
\begin{lemma} \label{topindlemma}
Let $ H \subset G $ be as above and let $ A $ be an $ H $-$ C^* $-algebra. Then $ \algind_H^G(A) $ 
is contained in $ C_b(G, A) $. 
\end{lemma} 
\proof Clearly it suffices to show $ p_\alpha C(G, A)^H \subset C_b(G, A) $ for $ \alpha \in \Irr(G)/\Irr(H) $.
That is, for $ f = (f_t)_{t \in \alpha} $ in $ C(G,A)^H $ with $ f_t \in C(G)_t \otimes A $ we have 
to show that $ ||f_t||_{t \in \alpha} $ is bounded. Upon replacing $ f $ by $ f^*f $ we may assume that 
all $ f_t $ are positive. Consider the element
$$ 
F = (\id \otimes \alpha)(f) = (\rho \otimes \id)(f) \in \prod_{r \in \alpha, s \in \Irr(H)} C(G)_r \otimes C(G)_s \otimes A
$$
Let us write $ || F ||_{r,s} $ for the norm of the restriction of $ F $ to $ C(G)_r \otimes C(G)_s \otimes A $. 
Since $ \alpha $ is isometric we have $ ||F||_{r,s} \leq ||F||_{r, \epsilon} = ||f_r|| $ for all $ s \in \Irr(H) $, where 
$ \epsilon $ denotes the trivial corepresentation. Next observe that we have $ ||(\rho \otimes \id)(f_t) ||_{r,s} = ||f_t|| $ 
provided $ t \subset r \otimes s $. From $ F = (\rho \otimes \id)(f) $ and the fact that 
all $ f_r $ are positive we obtain the estimate $ ||f_t|| \leq ||F||_{r,s} $ for all $ r \in \alpha, s \in \Irr(H) $ such that 
$ t \subset r \otimes s $. We conclude $ ||f_t|| \leq ||F||_{r,s} \leq ||F||_{r, \epsilon} = ||f_r|| $. 
This shows that $ ||f_t|| $ is in fact independent of $ t $, and in particular the family $ ||f_t||_{t \in \alpha} $ is bounded. \qed \\
Lemma \ref{topindlemma} shows that the canonical projection $ p_\alpha \algind_H^G(A) \rightarrow C(G)_t \otimes A $ is injective for all 
$ t \in \alpha $. Moreover we see that $ p_\alpha \algind_H^G(A) $ is a $ C^* $-subalgebra of $ C_b(G,A) $ for all $ \alpha \in \Irr(G)/\Irr(H) $. 
\begin{definition} \label{definduced}
Let $ H \subset G $ be a quantum subgroup of the discrete quantum group $ G $ and let $ A $ be an $ H $-$ C^* $-algebra.  
The induced $ C^* $-algebra $ \ind_H^G(A) $ is the closure of $ \algind_H^G(A) $ inside $ C_b(G, A) \subset M(C_0(G) \otimes A) $. 
\end{definition} 
Our above considerations imply that $ \ind_H^G(A) $ is simply the $ c_0 $-direct sum of the algebras $ p_\alpha \algind_H^G(A) $. 
Moreover from the fact that $ \algind_H^G(A) $ is a $ G $-algebra we see easily that $ \ind_H^G(A) $ is a $ G $-$ C^* $-algebra 
in a canonical way. \\
In the case of $ A = \mathbb{C} $ with the trivial action we write $ \ind_H^G(\mathbb{C}) = C_0(G/H) $ for the induced 
$ C^* $-algebra. Observe in particular that $ C_0(G/H) $ is a direct sum of finite dimensional $ C^* $-algebras. \\
Let us briefly compare definition \ref{definduced} with the general construction of induced $ C^* $-algebras by Vaes \cite{Vaesimprimitivity}. 
\begin{prop}
Let $ G $ be a discrete quantum group and let $ H \subset G $ be a quantum subgroup. If $ A $ is an $ H $-$ C^* $-algebra then 
the $ C^*$-algebra $ \ind_H^G(A) $ defined above is $ G $-equivariantly isomorphic to the induced $ C^* $-algebra defined by Vaes.
\end{prop}
\proof It suffices to check that the $ C^* $-algebra $ \ind_H^G(A) $ satisfies the conditions stated in theorem 7.2 of \cite{Vaesimprimitivity}. \\
We have already remarked above that $ \Delta \otimes \id: \ind_H^G(A) \rightarrow M(C_0(G) \otimes \ind_H^G(A)) $ yields 
a well-defined continuous coaction. Consider the algebra
$$
\tilde{A} = \{X \in M(\KH(l^2(G)) \otimes A) \mid X \in (L^\infty(G)' \otimes 1)' \;\text{and}\; (\rho \otimes \id)(X) = (\id \otimes \alpha)(X) \}
$$
used in \cite{Vaesimprimitivity}. For $ \alpha \in \Irr(G)/\Irr(H) $ and $ a \in A $ we have that 
$ X(p_\alpha \otimes a) $ and $ (p_\alpha \otimes a)X $ are contained in $ C_b(G, A) $ for all $ X \in \tilde{A} $. It follows that 
$ \algind_H^G(A) \subset \tilde{A} $ is strong* dense. We conclude that $ \Delta \otimes \id $ extends to 
a well-defined $ * $-homomorphism $ \tilde{A} \rightarrow M(C_0(G) \otimes \ind_H^G(A)) $ which is strictly continuous on 
the unit ball of $ \tilde{A} $. \\
Finally note that $ [\alpha(A)(l^2(H) \otimes A)] = l^2(H) \otimes A $ because $ \alpha $ is a coaction. 
Moreover we have $ [\CH[G] \cdot l^2(H)] = l^2(G) $ where we consider the left action of $ \mathbb{C}[G] $ 
corresponding to the regular coaction on $ l^2(G) $. Since we may view $ \alpha(A) $ as a subalgebra of 
$ \ind_H^G(A) $ we obtain
\begin{align*}
[l^2(G) \otimes A] &= [\CH[G] \cdot (l^2(H) \otimes A)] = [\CH[G] \cdot (\alpha(A)(l^2(H) \otimes A))] \\
&\subset [\CH[G] \cdot (\ind_H^G(A)(l^2(H) \otimes A))] 
\subset [\ind_H^G(A)(l^2(G) \otimes A)] 
\end{align*}
so that $ \ind_H^G(A) \subset \tilde{A} $ is nondegenerate. \qed \\ 
Next we shall analyze the Hilbert space associated to the homogeneous space $ G/H $. Since the Haar state for $ \mathbb{C}[G] $ restricts to 
the Haar state for $ \mathbb{C}[H] $ we have a canonical inclusion $ l^2(H) \subset l^2(G) $.  This yields a conditional expectation 
$ E: C^*_\red(G) \rightarrow C^*_\red(H) $ such that $ \hat{\phi} E = \hat{\phi} $ where $ \hat{\phi} $ denotes 
the Haar state of $ \CH[G] $. We have $ E(\CH[G]_r) = 0 $ if $ r \notin \Irr(H) $, and it is easily checked that 
$$
(E \otimes \id)\Delta = \Delta E = (\id\otimes E)\Delta. 
$$ 
Using $ E $ we define $ l^2(G/H) $ as the separated completion of $ \CH[G] $ with respect to the inner product 
$$
\bra x, y \ket = \epsilon E(x^*y). 
$$
The canonical left action of $ \CH[G] $ extends to a unital $ * $-homomorphism $ C^*_\max(G) \rightarrow \LH(l^2(G/H))$. We denote 
by $ l^2(G/H)_\alpha $ the image of $ \CH[G]_\alpha $ in $ l^2(G/H) $. \\
We want to show that the canonical map $ \CH[G/H] \rightarrow l^2(G/H) $ is injective. For this we will do computations in $ \CH[G] $ 
using matrix coefficients of corepresentations, and we will identify $ C(G)_r \cong \LH(\H_r) $ where $ \H_r $ is the underlying Hilbert 
space of $ r \in \Irr(G) $. There is a unique 
$ \overline{r} \in \Irr(G) $ such that $ r \otimes \overline{r} $ contains the trivial corepresentation $ \epsilon $. 
We denote by $ t_r: \epsilon \rightarrow r \otimes \overline{r} $ a morphism such that $ t_r^* t_r = \dim_q(r) 1 $ 
where $ \dim_q(r) $ is the quantum dimension of $ r $. Such a morphism is unique up to a phase, and 
there is a corresponding morphism $ s_r: \epsilon \rightarrow \overline{r} \otimes r $ such that 
$ (\id \otimes s_r^*)(t_r\otimes \id) = \id $. We also have $ (s_r^* \otimes \id)(\id \otimes t_r) = \id $ and 
$ s_r^* s_r = \dim_q(r) 1 $. That is, we can take $ t_{\overline{r}} = s_r $ and $ s_{\overline{r}} = t_r $. 
For $ x \in \CH[G]_{\overline{r}} $ and $ y \in \CH[G]_r $ the Schur orthogonality relations 
become
$$
\hat{\phi}(x y) = \frac{1}{\dim_q(r)} \bra x \otimes y, s_r s_r^* \ket
$$ 
in our notation. 
\begin{lemma} \label{GNSfaithful}
The canonical map $ \CH[G/H] \rightarrow l^2(G/H)$ is injective. 
\end{lemma}
\proof If $ x \in \CH[G]_r $ and $y \in \CH[G]_s $ with $ r, s \in \Irr(G) $ in different classes modulo $ \Irr(H) $, then 
$ x^*y \in \CH[G]_{\overline{r}} \cdot \CH[G]_s $ decomposes inside coefficients spaces $ \CH[G]_t $ with $ t \notin \Irr(H) $. 
This implies that the subspaces $ l^2(G/H)_\alpha $ for $ \alpha \in \Irr(G)/\Irr(H)$ are pairwise orthogonal. \\ 
In particular it suffices to check injectivity for fixed $ \alpha \in \Irr(G)/\Irr(H) $. We let $ r \in \alpha $ and recall 
that $ \CH[G/H]_\alpha = \tau(\CH[G]_r) $ is spanned by the elements from $ \CH[G]_r $. Let $ y $ be an element 
of $ \CH[G]_r $ such that $ \epsilon E(xy) = 0 $ for every $ x \in \CH[G]_{\overline{r}} $. We have to show $ y = 0 $
in $ \CH[G/H] $. \\
Let us choose a family of isometric morphisms $ T_i: s_i \rightarrow \overline{r} \otimes r $ for $ i = 1, \dots, N $ such 
that $ T_j^* T_i = \delta_{ij} $ and 
$$ 
q = \sum_{i = 1}^N  T_i T_i^* \in \LH(\H_{\overline{r}} \otimes \H_r) 
$$ 
is maximal such that $ s_i \in \Irr(H) $ for all $ i $. This allows us to compute 
\begin{align*}
\epsilon E(xy) = \bra E(xy), 1 \ket = \bra xy, p \ket = \sum_{i = 1}^N \bra x \otimes y, T_i T_i^* \ket = \bra x \otimes y, q \ket = \bra x, g \ket,
\end{align*}
where $ p \in C(G) $ is the orthogonal projection corresponding to $ \Irr(H) \subset \Irr(G) $ and 
$ g = (\id \otimes y)(q) \in \LH(\H_{\overline{r}}) $. \\
For each $ T_i $ we consider the morphism $ R_i = (\id \otimes T_i^*)(t_r \otimes \id):r \rightarrow r \otimes s_i $, and we define 
$ z \in \CH[G] \odot \CH[H] $ by 
\begin{equation*}
\bra z, f \otimes h \ket = \sum_{i = 1}^N \bra y, R_i^* (f \otimes h) R_i \ket.
\end{equation*}
If $ \mu $ denotes the multiplication of $ \CH[G]$, then $ \mu(z) $ is obtained by decomposing $ r \otimes s_i $ into irreducibles, and 
since $ r $ is irreducible there exists $ \lambda > 0 $ such that $ \lambda \mu(z) = y $. If we write 
$ y' = (\id \otimes \epsilon)(z) \in \LH(\H_r)^* = \CH[G]_r $ this implies $ y = \lambda y' $ in $ \CH[G/H] $.
To obtain an expression for $ y' $ we compute 
\begin{align*}
\sum_{i = 1}^N R_i^*(f \otimes 1) R_i &= \sum_{i = 1}^N (t_r^* \otimes \id)(\id \otimes T_i)(f \otimes 1)(\id \otimes T_i^*)(t_r \otimes \id) \\ 
&= (t_r^* \otimes \id)(f \otimes q)(t_r \otimes \id),
\end{align*}
so that $ \bra y', f \ket = t_r^*(f \otimes g)t_r $. Now we compute
\begin{align*}
\dim_q(r) \hat{\phi}(xy') &= \bra x \otimes y', s_r s_r^* \ket \\
&= \bra x, (\id\otimes t_r^*)(s_r s_r^* \otimes g)(\id\otimes t_r) \ket \\ 
&= \bra x, g \ket = \epsilon E(xy) = 0.
\end{align*}
for all $ x \in \CH[G]_{\overline{r}} $. Since $ \hat{\phi} $ is faithful on $ \CH[G] $ this yields $ y' = 0 $, and we 
conclude $ y = 0 $ in $ \CH[G/H] $. \qed \\
From lemma \ref{GNSfaithful} it follows in particular that the elements of $ C_c(G/H) \subset C(G) $ correspond to linear 
functionals on $ \CH[G] $ of the form $ d_E(x) $ where 
$$ 
d_H(x)(y) = \epsilon E(S(x)y). 
$$
Indeed, an element $ f \in C_c(G/H) $ corresponds to a linear form on $ \CH[G/H] $ which is supported on a finite subset of 
$ \Irr(G)/\Irr(H) $. Hence $ f $ can be written as the scalar product by a vector in $ l^2(G/H) $. Using lemma \ref{GNSfaithful} and the 
fact that the spaces $ \CH[G]_\alpha $ for $ \alpha \in \Irr(G)/\Irr(H) $ are finite dimensional we obtain the claim. \\
In this way we obtain a linear isomorphism $ d_H: \CH[G/H] \rightarrow C_c(G/H) $. This allows us to define an inner product on $ C_c(G/H) $ 
such that $ d_H $ is unitary. In the case of the trivial subgroup we write 
$$ 
d(x)(y) = \hat{\phi}(S(x)y),
$$ 
and note that $ d: \CH[G] \rightarrow C_c(G) $ is unitary with respect to the standard scalar products 
$$ 
\bra x, y \ket = \hat{\phi}(x^* y), \qquad \bra f, g \ket = \phi(f^*g)
$$
on $ \CH[G] $ and $ C_c(G) $, respectively. \\
It is easy to check that the comultiplication of $ \CH[G] $ induces a linear map $ \Delta: \CH[G/H] \rightarrow \CH[G/H] \odot \CH[G/H] $. 
Using the identification $ C(G/H) \cong \CH[G/H]^* $ we define the regular action of $ C(G/H) $ on $ \CH[G/H] $ by 
$$ 
\lambda(f)(x) = \bra f, x_{(1)} \ket x_{(2)}  
$$ 
for $ f \in C_c(G) $ and $ x \in \CH[G/H] $. This action turns $ \CH[G/H] $ into a left $ C(G/H) $-module. 
\begin{lemma}\label{lem_homogen_rep}
The regular action of $ C(G/H) $ on $ \CH[G/H] $ induces a faithful $ * $-representation $ \lambda: C_b(G/H) \rightarrow \LH(l^2(G/H)) $. 
Under this representation all elements of $ C_0(G/H) $ are mapped into $ \KH(l^2(G/H)) $.
\end{lemma}
\proof Under the isomorphism $ d: \CH[G] \rightarrow C_c(G) $ the canonical left action of an element $ f \in C(G) $ on $ C_c(G) $ is 
identified with the action $ \lambda(f)(x) = \bra f, x_{(1)} \ket x_{(2)} $ on $ \CH[G] $, and we have
$$ 
\hat{\phi}(\lambda(f^*)(x)^* y) = \hat{\phi}(x^* \lambda(f)(y)) 
$$ 
for all $ x, y \in \CH[G] $. If $ f \in C(G/H) \subset C(G) $ then the relation 
\begin{align*}
\lambda(f)(y) z &= \bra f, y_{(1)} \ket y_{(2)} z = \bra f, y_{(1)} z_{(1)} \ket y_{(2)} z_{(2)} = \lambda(f)(yz)
\end{align*}
for $ y \in \CH[G] $ and $ z \in \CH[H] $ shows 
$$ 
\hat{\phi}((\lambda(f^*)(x))^* y z) = \hat{\phi}(x^* \lambda(f)(y) z).  
$$ 
Since $ \hat{\phi} $ is faithful on $ \CH[H] $ we conclude  
$$
E(x^* \lambda(f)(y)) = E(\lambda(f^*)(x)^* y)
$$ 
for $ x, y \in \CH[G] $, and this implies 
$$ 
\epsilon E(x^* \lambda(f)(y)) = \epsilon E(\lambda(f^*)(x)^* y)
$$
for all $ x,y \in \CH[G/H] $. 
We conclude in particular that $ \lambda $ induces $ * $-homomorphisms $ p_\alpha C(G/H) \rightarrow \LH(l^2(G/H)_\alpha) $ 
for all $ \alpha \in \Irr(G)/\Irr(H) $, and taking the direct sum of these $ * $-homomorphisms yields the desired $ * $-representation 
$ \lambda: C_b(G/H) \rightarrow \LH(l^2(G/H)) $. Moreover, any element of $ p_\alpha C(G/H) $ for $ \alpha \in \Irr(G)/\Irr(H) $ acts as 
a finite rank operator because $ l^2(G/H)_\alpha = \CH[G/H]_\alpha $ is finite dimensional. It follows that $ \lambda $ maps the elements of 
$ C_0(G/H) $ to compact operators on $ l^2(G/H) $. \\
To prove faithfulness, observe that $ \lambda(f) = 0 \in \LH(l^2(G/H)) $ implies 
$$ 
\bra f, x \ket = \bra f, x_{(1)} \ket \epsilon(x_{(2)}) = \epsilon(\lambda(f)(x)) = 0
$$ 
for all $ x \in \CH[G/H] \subset l^2(G/H) $. Since $ C(G/H) $ identifies with the dual space of $ \CH[G/H] $ we conclude $ f = 0 $. \qed \\
Let us consider the left action of $ C_c(G) $ on $ C_c(G/H) $ given by 
$$
\mu(f)(h) = f_{(1)} h S(f_{(2)})
$$
for $ f \in C_c(G) $ and $ h \in C_c(G/H) $. 
Under the isomorphism $ d_E: \CH[G/H] \rightarrow C_c(G/H) $ defined above, this action identifies with the left action of 
$ C_c(G) $ on $ \CH[G/H] $ given by 
$$
\nu(f)(x) = \bra f_{(2)}, S(x_{(3)}) \ket x_{(2)} \bra f_{(1)}, x_{(1)} \ket.  
$$
It is straightforward to check that this induces a nondegenerate $ * $-homomorphism $ \nu: C_0(G) \rightarrow \LH(l^2(G/H)) $. \\
The $ G $-$ C^* $-algebra $ C_0(G/H) $ becomes a $ \DD(G) $-$ C^* $-algebra with the adjoint coaction 
$ C_0(G/H) \rightarrow M(C^*_\red(G) \otimes C_0(G/H)) $ given by 
$$
\gamma(f) = \hat{W}^*(1 \otimes f) \hat{W}
$$
where we view $ C_0(G/H) $ as a $ C^* $-subalgebra of $ C_b(G) \subset \LH(l^2(G)) $. 
\begin{lemma} 
The left regular action of $ C^*_\max(G) $ and the action $ \nu: C_0(G) \rightarrow \LH(l^2(G/H)) $ turn $ l^2(G/H) $ into a 
$ \DD(G) $-Hilbert space such that the $ * $-representation $ \lambda: C_0(G/H) \rightarrow \LH(l^2(G/H)) $ is covariant. 
\end{lemma} 
\proof The coaction $ \delta: l^2(G/H) \rightarrow M(C_0(G) \otimes l^2(G/H)) $ corresponding to the left action of 
$ C^*_\max(G) $ on $ l^2(G/H) $ is given by $ \Delta $ on $ C_c(G/H) \subset l^2(G/H) $. Using this fact it is straightforward to 
check the Yetter-Drinfeld compatibility condition and that 
$ \lambda: C_0(G/H) \rightarrow \LH(l^2(G/H)) $ is covariant. \qed \\
Let $ \H $ be a Hilbert space with scalar product $ \bra \;, \; \ket $ and let $ \S(\H) \subset \H $ be a dense linear subspace. We may 
view the scalar product on $ \S(\H) $ as a bilinear map $ \bra \;, \; \ket: \overline{\S(\H)} \times \S(\H) \rightarrow \mathbb{C} $ 
where $ \overline{\S(\H)} $ is the conjugate of $ \S(\H) $. That is, $ \S(\H) = \overline{\S(\H)} $ is the vector 
space equipped with the same addition but with the complex conjugate scalar multiplication. 
We let $ \mathbb{F}(\H) = \S(\H) \odot \overline{\S(\H)} $ be the $ * $-algebra with multiplication 
$$
(x_1 \odot y_1)(x_2 \odot y_2) = x_1 \bra y_1, x_2 \ket \odot y_2
$$
and $ * $-operation
$$
(x \odot y)^* = y \odot x. 
$$
Note that, by slight abuse of notation, $ \mathbb{F}(\H) $ depends on $ \S(\H) $ and not only $ \H $. \\
We may view $ \mathbb{F}(\H) $ as a space of kernels for certain compact operators on $ \H $, 
more precisely, there is an injective $ * $-homomorphism $ \iota: \mathbb{F}(\H) \rightarrow \KH(\H) $ given by 
$$
\iota(x \odot y)(\xi) = x \bra y, \xi \ket.  
$$
We apply these constructions to the homogeneous space associated to a discrete quantum group $ G $ 
and a quantum subgroup $ H \subset G $. In the Hilbert space $ l^2(G/H) $ we have the dense linear subspace $ \CH[G/H] $.
Let $ \gamma: \CH[G/H] \rightarrow \overline{\CH[G/H]} $ be the linear map given by $ \gamma(x) = \overline{S(x)^*} $. 
It is easy to check that $ \gamma $ is a well-defined linear isomorphism. \\
We may therefore identify the algebra $ \mathbb{F}(l^2(G/H)) $ with 
$ \mathbb{C}[G/H] \odot \mathbb{C}[G/H] $. The actions of $ \CH[G] $ and $ C_c(G) $ on $ \KH(l^2(G/H)) $ induced from the 
Yetter-Drinfeld structure of $ l^2(G/H) $ preserve $ \mathbb{F}(l^2(G/H)) $, and correspond to the diagonal actions on 
$ \mathbb{C}[G/H] \odot \mathbb{C}[G/H] $. \\
The operators in $ \mathbb{F}(l^2(G/H)) $ of the form $ \Delta(x) $ with $ x \in \CH[G/H] $ form a subalgebra, and 
$ \Delta(\CH[G/H]) \subset \mathbb{F}(l^2(G/H)) $ is closed under the the actions of $ \mathbb{C}[G] $ 
and $ C_c(G) $. It is easy to check that the closure of $ \Delta(\CH[G/H]) $ inside $ \KH(l^2(G/H)) $ identifies 
with $ C_0(G/H) $. \\
In a similar way we may consider operator kernels defining operators between $ l^2(G/K) $ and $ l^2(G/H) $ if $ H $ and $ K $ 
are two different quantum subgroups of $ G $. Later we will need in particular the following fact. 
\begin{lemma} \label{moritalemma}
Let $ G $ be a discrete quantum group and let $ H \subset G $ be a quantum subgroup. 
Consider the vector space 
$$ 
L = \Delta(\mathbb{C}[G])(\mathbb{C}[H] \odot \mathbb{C}[H]) 
$$
inside $ \mathbb{C}[G] \odot \mathbb{C}[G] $.  
If $ \tau: \CH[G] \rightarrow \CH[G/H] $ denotes the canonical projection, then the closure $ K $ of 
the $ * $-algebra
$$ 
\begin{pmatrix} 
L & (\id \odot \tau)(L) \\
(\tau \odot \id)(L) & (\tau \odot \tau)(L)
\end{pmatrix}
$$
inside $ \KH(l^2(G) \oplus l^2(G/H)) $ is $ \DD(G) $-equivariantly Morita equivalent to $ C_0(G/H) $.
\end{lemma} 
\proof It is easy to check that $ K $ is indeed closed under multiplication and taking adjoints. 
Let $ p \in M(K) $ be the projection onto $ l^2(G/H) $. This projection is invariant under the action of the 
Drinfeld double in the sense that $ \gamma(p) = 1 \otimes p $ if $ \gamma: K \rightarrow M(C^\red_0(\DD(G)) \otimes K) $ denotes the coaction 
of $ C_0^\red(\DD(G)) $ induced from $ l^2(G) \oplus l^2(G/H) $. Moreover, by our above remarks we know that $ p K p $ is isomorphic to 
$ C_0(G/H) $. Hence the claim follows from the fact that $ p K p $ is a full corner in $ K $, that is $ [K p K] = K $. \qed

\section{Divisible quantum subgroups} \label{secdivisible}

In this section we discuss a certain class of quantum subgroups of discrete quantum groups which we call divisible. We will later make use 
of divisibility in some arguments. \\
Let $ G $ be a discrete quantum group and let $ H \subset G $ be a quantum subgroup. In the same way as we defined the $ C^* $-algebra 
$ C_0(G/H) $ in section \ref{secquantumsub}, we may construct the homogeneous space $ C_0(H \backslash G) $ corresponding to left translations 
of $ H $ on $ G $. That is, on the algebraic level we set
\begin{equation*}
C(H \backslash G) = \{f \in C(G) \mid (\pi \odot \id)\Delta(f) = 1 \odot f \}
\end{equation*}
where $ \pi: C(G) \rightarrow C(H) $ is the canonical projection map, and $ C_0(H \backslash G) $ is the completion of 
$$ 
C_c(H \backslash G) = \bigoplus_{\alpha \in \Irr(H) \backslash \Irr(G)} p_\alpha C(H \backslash G) 
$$
inside $ \LH(l^2(G)) $. Here $ \Irr(H) \backslash \Irr(G) $ is the quotient of $ \Irr(G) $ with respect to the equivalence relation given by 
$ r \sim s $ iff $ s \otimes \overline{r} $ contains an element $ t \in \Irr(H) $, and $ p_\alpha \in C_b(G) $ denotes 
the projection determined by the coset $ \alpha $. 
\begin{definition} \label{defdivisible}
Let $ G $ be a discrete quantum group and let $ H \subset G $ be a quantum subgroup. We say that $ H \subset G $ is divisible if there exists 
an $ H $-equivariant $ * $-isomorphism 
$$ 
C_0(G) \cong C_0(H) \otimes C_0(H \backslash G) 
$$ 
with respect to the restricted coaction on the left hand side and the coaction given by comultiplication on $ C_0(H) $ on the right hand side. 
\end{definition} 
Roughly speaking, the property of divisibility corresponds to the existence of a section $ H \backslash G \rightarrow G $ of the canonical
quotient map. It is clear that every inclusion $ H \subset G $ of classical discrete groups is divisible. \\
Recall from section \ref{secquantumsub} that the canonical projection $ \CH[G] \rightarrow \CH[G/H] $ is a coalgebra homomorphism 
in a natural way. We define
$$
\CH[H \backslash G] = \mathbb{C} \odot_{\CH[H]} \CH[G] 
$$
and note that $ \CH[H \backslash G]^* $ identifies with $ C(H \backslash G) $. The coalgebra $ \CH[G] $ is equipped with the 
antilinear involution $ \rho(x) = S^{-1}(x)^* $. This formula also determines involutions on $ \CH[G/H] $ and $ \CH[H \backslash G] $. We 
will say that a coalgebra homomorphism is involutive if it commutes with these involutions. \\
The following lemma shows that divisibility is quite a strong property in the quantum case, and that it is symmetric 
with respect to taking left or right quotients. 
\begin{lemma} \label{chardivisible}
Let $ G $ be a discrete quantum group and let $ H \subset G $ be a quantum subgroup. Then the following conditions are equivalent. 
\begin{bnum}
\item[a)] $ H \subset G $ is divisible. 
\item[b)] There exists an involutive isomorphism $ \CH[G] \cong \CH[H] \odot \CH[G \backslash H] $ of coalgebras and left $ \CH[H] $-modules. 
\item[c)] For each $ \alpha \in \Irr(H) \backslash \Irr(G) $ there exists a corepresentation $ l = l(\alpha) \in \alpha $ such that 
$ s \otimes l $ is irreducible for all $ s \in \Irr(H) $.
\item[d)] For each $ \alpha \in \Irr(G)/\Irr(H)$ there exists a corepresentation $ r = r(\alpha) \in \alpha $ such that 
$ r \otimes s $ is irreducible for all $ s \in \Irr(H) $.
\item[e)] There exists an involutive isomorphism $ \CH[G] \cong \CH[G/H] \odot \CH[H] $ of coalgebras and right $ \CH[H] $-modules. 
\item[f)] There exists a right $ C_0(H) $-colinear $ * $-isomorphism 
$$ 
C_0(G) \cong C_0(G/H) \otimes C_0(H) 
$$ 
with respect to the restricted coaction on the left hand side and the coaction given by comultiplication on $ C_0(H) $ on the 
right hand side. 
\end{bnum}
\end{lemma}
\proof $ a) \Rightarrow b) $ Let $ \sigma: C_0(G) \rightarrow C_0(H) \otimes C_0(H \backslash G) $ be an $ H $-equivariant 
$ * $-isomorphism. Then every matrix block in $ C_0(G) $ is mapped to finitely many matrix blocks inside 
$ C_0(H) \otimes C_0(H \backslash G) $. It follows that there exists an isomorphism $ \sigma_*: \CH[H] \odot \CH[H \backslash G] \rightarrow 
\CH[G] $ which induces $ \sigma $ by dualizing in the sense that 
$$
\bra \sigma_*(x \otimes y), f \ket = \bra x \otimes y, \sigma(f) \ket 
$$ 
for all $ x \in \CH[H], y \in \CH[H \backslash G] $ and $ f \in C_0(G) $. Using nondegeneracy of the canonical pairings
we see that $ \sigma_* $ is a left $ \CH[H] $-linear involutive coalgebra homomorphism. \\
$ b) \Rightarrow a) $ If $ \sigma: \CH[H] \odot \CH[H \backslash G] \rightarrow 
\CH[G] $ is a $ \CH[H] $-linear involutive coalgebra isomorphism, then dualizing yields a $ \CH[H] $-linear 
unital $ * $-isomorphism $ \hat{\sigma}: C(G) = \M(C_c(G)) \rightarrow \M(C_c(H) \odot C_c(H \backslash G)) $. 
Since both sides are direct products of matrix algebras the map $ \hat{\sigma} $ restricts to a $ C_0(H) $-colinear 
$ * $-isomorphism $ C_0(G) \cong C_0(H) \otimes C_0(H \backslash G) $. \\
$ b) \Rightarrow c) $ Assume we are given an involutive $ \CH[H] $-linear coalgebra isomorphism $ \CH[G] \cong \CH[H] \odot \CH[H \backslash G] $.  
Let $ \alpha \in \Irr(H) \backslash \Irr(G) $ and $ l_1, \dots, l_k \in \Irr(G) $ 
be the corepresentations corresponding to the simple coalgebras in $ \CH 1 \odot \CH[H \backslash G]_\alpha $. Then the simple 
coalgebras in $ \CH[G]_\alpha $ are of the form $ \CH[G]_{s \otimes l_i} $ with $ s \in \Irr(H) $. Since $ \alpha $ is 
generated by a single corepresentation we must have $ k = 1 $ and may set $ l(\alpha) = l_1 \in \Irr(G) $. \\
$ c) \Rightarrow b) $ If $ l(\alpha) \in \Irr(G) $ is a representative of $ \alpha \in \Irr(H) \backslash \Irr(G) $, then the simple 
corepresentations of $ \CH[G] $ are of the form $ s \otimes l(\alpha) $ with $ s \in \Irr(H) $. The multiplication map $ \mu $ induces 
a left $ \CH[H] $-linear isomorphism
\begin{align*}
\bigoplus_{\alpha \in \Irr(H) \backslash \Irr(G)} \mathbb{C}[H] \odot \mathbb{C}[G]_{l(\alpha)} \cong \mathbb{C}[G]
\end{align*}
of coalgebras. Using $ \rho(xy) = S(xy)^* = S(x)^* S(y)^* = \rho(x) \rho(y) $ we see that $ \mu $ is involutive. The claim follows by 
observing that the canonical map 
$$
\bigoplus_{\alpha \in \Irr(H) \backslash \Irr(G)} \mathbb{C}[G]_{l(\alpha)} \rightarrow \CH[H \backslash G]
$$
is an involutive isomorphism of coalgebras. \\
$ c) \Leftrightarrow d) $ The conjugation of corepresentations induces a canonical bijection between $ \Irr(H) \backslash \Irr(G) $ and 
$ \Irr(G)/\Irr(H) $. If $ l = l(\alpha) $ is as in $ c) $ we may set $ r(\overline{\alpha}) = \overline{l} $ and vice versa. \\
Finally, the proof of $ d) \Leftrightarrow e) \Leftrightarrow f) $ is obtained in the same way as the equivalences
$ a) \Leftrightarrow b) \Leftrightarrow c) $ above. \qed \\
In order to illustrate definition \ref{defdivisible} let us discuss some examples of quantum subgroups and homogeneous spaces. 
As a first example, let us consider the even part $ H = \FH O_{ev}(Q) $ of the free orthogonal quantum group $ G = \FH O(Q) $ with 
$ Q \overline{Q} = \pm1 $. By definition, $ \CH[H] $ is the Hopf $ * $-subalgebra generated by products $ u_{ij} u_{kl} $ of two 
generators of $ \CH[G] $. In other words, if 
$$ 
\Irr(G) = \{u_k \mid k \in \mathbb{N}_0 \} 
$$ 
is the usual parametrization of the irreducible corepresentations of $ \FH O(Q) $, then $ \Irr(H) $ corresponds to corepresentations 
with even index. Note that one recovers as particular case the canonical morphisms $ SU_q(2) \rightarrow SO_q(3) $ when 
$ Q \in GL_2(\mathbb{C}) $. \\
Recall that the fusion rules of $ \FH O(Q)$ are given by 
\begin{align*}
u_k \otimes u_l \cong u_{|k - l|} \oplus u_{|k - l| + 2} \oplus \cdots \oplus u_{k + l - 2} \oplus u_{k + l}, \qquad 
\overline{u_k} \cong u_k.
\end{align*}
Hence in this case $ \Irr(G)/\Irr(H) $ consists of the two classes $ [1] = [u_0] $ and $ [u] = [u_1] $. Clearly 
$ \CH[G/H]_{[1]} \cong \CH $ is generated by $ 1 \in \CH[G] $. On the other hand, for $ r, s \in [u] $ we have 
$ \overline{r} \otimes s \in \Corep(H) $ so that $ \bra x, y \ket = \epsilon(x^* y) $ for $ x, y \in \CH[G/H]_{[u]} $. As a 
result we obtain $ \CH[G/H]_{[u]} \cong \CH u_{11} $, with the coalgebra structure which makes $ u_{11} $ 
a group-like element. In other words, we have 
$$ 
C_0(G/H) = C_c(G/H) \cong \mathbb{C} \oplus \mathbb{C}, 
$$
and it follows immediately that $ H $ is not a divisible subgroup. \\
We consider next the case of a free product $ G = G_0 * G_1 $ of two arbitrary discrete quantum groups $ G_0, G_1 $. It is known from 
\cite{Wangfree} that the irreducible corepresentations of $ G $ are obtained as alternating tensor products of nontrivial irreducible corepresentations of $ G_0 $ and $ G_1 $. We identify $ \Irr(G) $ with the set of alternating words $ \Irr(G_0) * \Irr(G_1) $. The fusion 
rules of $ G $ can be recursively derived from the ones of $ G_0 $ and $ G_1 $ as follows. If the word $ v $ ends in $ \Irr(G_0) $ and 
$ w \in \Irr(G) $ starts in $ \Irr(G_1) $ then $ v \otimes w $ is irreducible and corresponds to the concatenation $ vw $. The same holds 
with the roles of $ G_0 $ and $ G_1 $ reversed. If on the other hand $ v = v' r $ and $ w = s w' $ with $ r, s \in \Irr(G_i) $, then
\begin{align*}
v \otimes w \cong \bigoplus_{\epsilon \neq t \subset r \otimes s} v' t w' \oplus \delta_{\overline{r}, s} (v' \otimes w'),
\end{align*}
where the sum runs over all irreducibles corepresentations $ t \subset r \otimes s $ with multiplicities, according to the fusion rules of 
$ G_i $. It follows in particular that each $ G_i \subset G $ is a divisible quantum subgroup, with the corepresentations 
$ r(\alpha) $ for $ \alpha \in \Irr(G)/\Irr(G_i) $ corresponding to words not ending in $ G_i $, including the empty word. \\ 
We remark that one can construct a naive Bass-Serre tree $ Y $ from these fusion rules. The set of edges is 
the disjoint union $ Y^{(0)} = \Irr(G)/\Irr(G_0) \cup \Irr(G)/\Irr(G_1) $, and the set of vertices is $ Y^{(1)} = \Irr(G) $. 
The natural quotient maps $ \tau_j : \Irr(G) \to \Irr(G)/\Irr(G_j) $ are then interpreted as source and target maps, thus 
defining an oriented graph $ Y $. It is easy to check that this graph is a tree, and it coincides with the Bass-Serre tree if 
$ G_0, G_1 $ are classical groups. In the quantum case, the free product $ G $ does not act on $ Y $ in a natural way, and 
we will work instead with a quantum Bass-Serre tree in section \ref{secdirac}. Still, the naive tree $ Y $ described above 
turns out to be useful for $ K $-theoretic computations. \\
Finally, let us consider free orthogonal and unitary quantum groups. It was shown by Banica \cite{Banicaunitary} that for any 
$ Q \in GL_n(\mathbb{C}) $ such that $ Q \overline{Q} = \pm 1 $ there is an injective morphism 
$ \iota: \mathbb{C}[\mathbb{F}U(Q)] \rightarrow \mathbb{C}[\mathbb{F}O(Q) * \mathbb{Z}] $ of Hopf-$ * $-algebras given by 
$ \iota(u_{ij}) = u_{ij} z $ where $ z \in \mathbb{C}[z, z^{-1}] = \CH[\mathbb{Z}] $ is the generator. That is, 
$ \mathbb{F}U(Q) \subset \mathbb{F}O(Q) * \mathbb{Z} $ is a quantum subgroup. 
\begin{prop} \label{fuinc}
Let $ Q \in GL_n(\mathbb{C}) $ such that $ Q \overline{Q} = \pm 1 $. Then the quantum subgroup 
$ \mathbb{F}U(Q) \subset \mathbb{F}O(Q) * \mathbb{Z} $ is divisible. 
\end{prop} 
\proof Let us write $ G = \mathbb{F}O(Q) * \mathbb{Z} $ and $ H = \mathbb{F}U(Q) $ and identify $ \Irr(H) \subset \Irr(G) $. 
Moreover denote by $ u_n $ for $ n \in \mathbb{N}_0 $ the irreducible corepresentations of $ \mathbb{F}O(Q)$, and choose a generator 
$ z $ of $ \mathbb{Z} $. By definition, $ \Irr(H) $ is generated inside $ \Irr(G) $ by $ u_1 z $ and $ z^{-1} u_1 $. 
We have to study how tensor products by the generators decompose in $ \Irr(G) $. \\
Let $ w z^k u_n z^l $ be an alternating word in $ \Irr(G) $ with $ w $ possibly empty, $ k, l \in \mathbb{Z} $ and 
$ n > 0 $. According to the fusion rules explained above we have 
\begin{align*}
w z^k u_n z^l \otimes u_1 z &= w z^k u_n z^l u_1 z \qquad \text{if} \; l \neq 0 \\ 
w z^k u_n \otimes u_1 z &= w z^k u_{n - 1}z \oplus w z^k u_{n + 1} z \qquad \text{if} \; n > 1 \\
w z^k u_1 \otimes u_1 z &= w z^{k + 1} \oplus w z^k u_2 z. 
\end{align*} 
Similarly we have 
\begin{align*}
w u_n z^k \otimes z^{-1} u_1 &= w u_n z^{k - 1} u_1 \qquad \text{if} \; k \neq 1 \\
w u_n z \otimes z^{-1} u_1 &= w u_{n - 1} \oplus w u_{n + 1} \qquad \text{if} \; n > 1 \\
w u_1 z \otimes z^{-1} u_1 &= w \oplus w u_2
\end{align*}
for $ k \in \mathbb{Z}, n > 0 $ and a possibly empty word $ w $. \\ 
Let us describe $ \Irr(H) $ inside $ \Irr(G) $. Indexing irreducible corepresentations of $ \mathbb{F}O(Q)* \mathbb{Z} $ by words 
in $ u_n $, $ z^k $ with $ n \in \mathbb{N}^*, k \in \mathbb{Z}^* $, as explained above, we claim that the nontrivial irreducible
corepresentations of $ \mathbb{F}U(Q) $ correspond to the set $ W $ of words of the form 
\begin{align*}
z^{[\epsilon_0]_-}u_{n_1}z^{\epsilon_1}u_{n_2}z^{\epsilon_2}\cdots
u_{n_p}z^{[\epsilon_p]_+}
\end{align*}
where $ p \geq 1 $, $ \epsilon_i = \pm 1 $ and $ n_i \in \mathbb{N}^*$ for all $ i $ are such that 
$ \epsilon_{i + 1} = -(-1)^{n_{i + 1}} \epsilon_i $ for all $ 0 \leq i \leq p - 1 $ where we use the notation 
\begin{align*}
[-1]_- = -1, \qquad [1]_- = 0, \qquad [-1]_+ = 0, \qquad [1]_+ = 1.
\end{align*}
Indeed, it is easy to check that $ u_1 z = z^{[1]_-} u_1 z^{[1]_+} $ and $ z^{-1} u_1 = z^{[-1]_-} u_1 z^{[-1]_+} \in W $. 
Moreover one verifies that $ W $ is stable under tensoring on the left or right by $ u_1 z $ or $ z^{-1} u_1 $ and taking subobjects 
according to the fusion rules above, and that all words of $ W $ can be obtained in this way. For the rest of proof we only need the 
inclusion $ \Irr(H) \subset W $. \\
Let $ \alpha \in \Irr(G)/ \Irr(H) $ be given. According to lemma \ref{chardivisible}, it suffices to find $ r \in \alpha $ such that 
$ r \otimes t $ is irreducible for all $ t \in \Irr(H) $. Let us pick an element $ r $ in $ \alpha $ of minimal length as a word in 
$ \mathbb{N} * \mathbb{Z} $. If $ r $ has length $ 0 $ we have $ r = \epsilon $ and $ \alpha = \Corep(H) $, in particular the irreducibility 
property is clear in this case. If $ r = u_n $ with $ n > 0 $, then tensoring by 
$ (uz \otimes z^{-1} u)^{\otimes m} = u^{\otimes 2m} \in \Corep(H) $ shows that we can assume $ r = u_1 $. Tensoring on the 
right by $ u_1 z $ shows that we can also take $ r = z $, and the irreducibility property is again satisfied. \\
Now assume that $ r $ is of length at least $ 2 $. If $ r $ ends by $ u_n $ with $ n > 0 $, then we can achieve $ n = 1 $ as above, and 
tensoring by $ u_1 z $ yields a word of strictly smaller length in $ \alpha $. Hence $ r $ must be of the form $ w u_n z^k $ with 
$ k \in \mathbb{Z}^*, n \in \mathbb{N}^* $. In fact $ k \neq 1 $, because otherwise we find again a strictly shorter word in 
$ \alpha $ by tensoring with $ z^{-1} u $. Now it is clear from the fusion rules and the description of $ \Irr(H) $ above 
that the tensor products $ r \otimes t $ with $ t \in \Irr(H) $ are irreducible. \qed

\section{The Dirac and dual Dirac elements} \label{secdirac}

In this section we define the Bass-Serre quantum tree and the Dirac and dual Dirac elements associated to a free product of discrete 
quantum groups. We shall essentially follow Kasparov and Skandalis \cite{KSbuildings}, in addition we take into account the natural 
Yetter-Drinfeld structures in the quantum setting. In the sequel we assume that the reader is familiar with equivariant $ KK $-theory 
for quantum groups, our notation will follow \cite{NVpoincare}, \cite{Voigtbcfo}. \\
Let $ G_0 $ and $ G_1 $ be discrete quantum groups and let $ G = G_0 * G_1 $ be their free product. We write 
$ \tau_j: \mathbb{C}[G] \rightarrow \mathbb{C}[G/G_j] $ for the canonical projections. On the level of Hilbert spaces, we define an 
associated tree for $ G $ by 
$$ 
l^2(X) = l^2(X^{(0)}) \oplus l^2(X^{(1)}) 
$$ 
where 
$$ 
l^2(X^{(0)}) = l^2(G/G_0) \oplus l^2(G/G_1), \qquad l^2(X^{(1)}) = l^2(G). 
$$
Heuristically, this Hilbert space can be viewed as the space of $ l^2 $-summable functions on the simplices of the quantum tree $ X $ 
associated to the free product $ G $. However, we point out that $ X $ itself has no meaning here, it is only the Hilbert space $ l^2(X) $ 
with further structure which will be used in the sequel. Observe that $ l^2(X) $ contains a natural dense linear subspace 
$ \S(l^2(X)) = \mathbb{C}[X] $ given by
$$ 
\mathbb{C}[X] = \mathbb{C}[X^{(0)}] \oplus \mathbb{C}[X^{(1)}] = 
\mathbb{C}[G/G_0] \oplus \mathbb{C}[G/G_1] \oplus \mathbb{C}[G].   
$$ 
We will work with compact operators on $ l^2(X) $ obtained from kernels in $ \mathbb{C}[X] \odot \mathbb{C}[X] $ as explained in 
section \ref{secquantumsub}. \\
Let 
$$ 
E = \{(t_0,t_1) \in \mathbb{R}^2| t_0 + t_1 = 1 \}
$$
be the one-dimensional affine space and let 
$$
\Delta^1 = \{(t_0,t_1) \in E| t_0 \geq 0, t_1 \geq 0 \} \subset E 
$$
be the standard $ 1 $-simplex. As in \cite{KSbuildings} we denote by $ q: E \rightarrow \Delta^1 $ the projection to the nearest point. 
For a subset $ f \subset \{0, 1\} $ define 
$$
F_f = \{(t_0, t_1) \in \Delta^1| t_j = 0 \;\text{for}\; j \in \{0,1\} \setminus f \}. 
$$
Note that $ F_{\{0,1\}} = \Delta^1 $ and $ F_\emptyset = \emptyset $. Moreover define the open set $ \Omega_f \subset E $ as the interior 
in $ E $ of $ q^{-1}(F_f) $. Observe that $ \Omega_{\{0,1\}} = E $ and $ \Omega_\emptyset = \emptyset $. For the one-element subset 
$ \{j\} \subset \{0,1\} $ we shall abbreviate $ \Omega_{\{j\}} = \Omega_j $. \\
Let $ \mathbb{C}l_1 $ be the one-dimensional Clifford algebra and write 
$$ 
C_\tau(U) = C_0(U) \otimes \mathbb{C}l_1
$$ 
if $ U \subset E $ is an open subset. Our first aim is to define a graded $ C^* $-algebra 
$$ 
\A_X \subset C_\tau(E) \cotimes \hat{\KH}(l^2(X)) = C_0(E, \hat{\KH}(l^2(X))) \cotimes \mathbb{C}l_1 
$$ 
where $ C_\tau(E) $ carries the grading induced from the Clifford algebra, $ \cotimes $ denotes the graded tensor product and 
$ \hat{\KH}(l^2(X)) $ is the algebra $ \KH(l^2(X)) $ of compact operators on $ l^2(X) $ with the grading obtained by viewing 
$ l^2(X) = l^2(X^{(0)}) \oplus l^2(X^{(1)}) $ as a graded Hilbert space. In the case of classical groups, the algebra $ \A_X $ is defined 
by specifying support conditions in terms of matrix elements for $ \KH(l^2(X)) $. In the quantum case we have to proceed more indirectly, we 
shall assemble $ \A_X $ by putting together the pieces corresponding to the different regions $ \Omega_f $ for $ f \subset \{0,1\} $. \\
Firstly, let us consider $ \Omega_{\{0,1\}} = E $. Here we take the linear subspace 
$$ 
C_\tau(E) \odot \Delta(\mathbb{C}[G]) \subset C_\tau(E) \odot \mathbb{F}(l^2(X)). 
$$
As explained in section \ref{secquantumsub}, the closure of this space inside $ C_\tau(E) \cotimes \hat{\KH}(l^2(X)) $ 
is isomorphic to $ C_\tau(E) \otimes C_0(G) $. Heuristically, this contribution correponds to the matrix elements for 
$ \sigma, \eta \in X^{(1)} $ such that $ \sigma = \eta $. \\
Consider next the regions $ \Omega_{\{j\}} = \Omega_j $ for $ j = 0,1 $. Heuristically, we have a contribution coming from all simplices 
which intersect in the vertices corresponding to the homogeneous space $ G/G_j $. For a pair of simplices in $ X^{(1)} $ this leads to the space 
$$ 
C_\tau(\Omega_j) \odot \Delta(\mathbb{C}[G])(\mathbb{C}[G_j] \odot \mathbb{C}[G_j]) 
$$ 
inside $ C_\tau(E) \odot \mathbb{F}(l^2(X^{(1)})) $. For a pair of one simplex in $ X^{(1)} $ and one simplex in $ X^{(0)} $ we have
$$ 
C_\tau(\Omega_j) \odot (\tau_j \odot \id)\Delta(\mathbb{C}[G]) 
$$ 
inside $ C_\tau(E, \FH(l^2(X^{(0)}), l^2(X^{(1)})) $ and symmetrically 
$$
C_\tau(\Omega_j) \odot (\id \odot \tau_j)\Delta(\mathbb{C}[G]) 
$$ 
inside $ C_\tau(E, \FH(l^2(X^{(1)}), l^2(X^{(0)})) $. Finally, there is a contribution coming from a pair of vertices in $ X^{(0)} $ which 
gives
$$ 
C_\tau(\Omega_j) \odot (\tau_j \odot \tau_j)\Delta(\mathbb{C}[G]) 
$$ 
inside $ C_0(E, \FH(l^2(X^{(0)})) $. \\
We let $ \A_X $ be the closure of the sum of the above subspaces inside $ C_\tau(E) \cotimes \hat{\KH}(l^2(X)) $. It is straightforward to 
check that $ \A_X $ becomes a $ \DD(G) $-$ C^* $-algebra in a natural way. 
\begin{definition}
Let $ G_0, G_1 $ be discrete quantum groups and let $ G = G_0 * G_1 $ their free product. The Dirac element 
$ D \in KK^{\DD(G)}(\A_X, \mathbb{C}) $ is the composition of the canonical inclusion $ \A_X \rightarrow C_\tau(E) \cotimes \hat{\KH}(l^2(X)) $ 
with the Bott periodicity isomorphism and the equivariant Morita equivalence $ \hat{\KH}(l^2(X)) \sim_M \mathbb{C} $. 
\end{definition}
Note that $ l^2(X) $ is a graded $ \DD(G) $-Hilbert space, so that $ \hat{\KH}(l^2(X)) $ is in fact $ \DD(G) $-equivariantly Morita equivalent 
to $ \mathbb{C} $. \\
The $ \DD(G) $-$ C^* $-algebra $ \A_X $ is $ KK^{\DD(G)} $-equivalent to an ungraded $ \DD(G) $-$ C^* $-algebra. More precisely, let 
$ \B_X \subset C_0(E) \otimes \KH(l^2(X)) $ be the ungraded $ C^* $-algebra defined by the same procedure as above, but replacing $ C_\tau(U) $ 
by $ C_0(U) $ in all steps. Then $ \P = C_0(\mathbb{R}) \otimes \B_X = \Sigma \B_X $ is isomorphic to $ \A_X $ in $ KK^{\DD(G)} $. Indeed, the 
grading on $ \hat{\KH}(l^2(X)) $ is even and hence $ C_\tau(E) \cotimes \hat{\KH}(l^2(X)) $ is isomorphic to 
$ C_0(E) \otimes \mathbb{C}l_1 \otimes \KH(l^2(X)) $, where $ \mathbb{C}l_1 \otimes \KH(l^2(X)) $ carries the standard odd grading, see 
corollary 14.5.3 in \cite{Blackadar}. This isomorphism identifies $ \A_X $ with $ \mathbb{C}l_1 \otimes \B_X $. Since $ \mathbb{C}l_1 $ is 
$ KK $-equivalent to $ C_0(\mathbb{R}) $ this yields the claim. \\
In order to analyze the structure of $ \P $ it is useful to consider the projection $ \P \rightarrow \Sigma C(\Delta^1, C_0(G)) $ obtained 
from restriction of functions in $ \B_X $ to $ \Delta^1 \subset E $. We otain an extension
$$
\xymatrix{
0 \ar@{->}[r] & I_0 \oplus I_1 \ar@{->}[r] & \P \ar@{->}[r] & \Sigma C(\Delta^1, C_0(G)) \ar@{->}[r] & 0 
}
$$
of $ \DD(G) $-$ C^* $-algebras, where $ I_j \subset \P $ is the ideal corresponding to the open sets $ \Omega_j $. One checks easily that 
this extension has a $ \DD(G) $-equivariant completely positive splitting. Moreover, it follows from lemma \ref{moritalemma} that $ I_j $ 
is $ KK^{\DD(G)} $-equivalent to $ C_0(G/G_j) $ for $ j = 0,1 $. \\
Next we shall define the $ \gamma $-element for $ G $, compare \cite{JVkam}, \cite{Vergniouxkam}. If $ H $ is a quantum group let us write 
$ \mathbb{C}[H \setminus \{e\}] $ for the subspace of $ \mathbb{C}[H] $ spanned by the matrix coefficients of all nontrivial corepresentations 
of $ H $. Similarly, we write $ l^2(H \setminus \{e\}) \subset l^2(H) $ for the orthogonal complement of $ \mathbb{C} 1 $. That is, we 
have a direct sum decomposition 
$$
l^2(H) = \mathbb{C} 1 \oplus l^2(H \setminus \{e\})
$$
of the Hilbert space $ l^2(H) $. With this notation, the underlying vector space of the Hopf algebra $ \mathbb{C}[G] $ for the free product 
$ G = G_0 * G_1 $ decomposes as
$$
\mathbb{C}[G] = \mathbb{C} 1 \oplus \bigoplus_{n = 1}^\infty \; \bigoplus_{(i_1, \dots, i_n) \in I_n} \mathbb{C}[G_{i_1}\setminus \{e\}] 
\odot \cdots \odot \mathbb{C}[G_{i_n}\setminus \{e\}]
$$
where 
$$
I_n = \{(i_1, \dots, i_n) \in \{0,1\}^n| i_k \neq i_{k + 1}\; \text{for all}\; k \}.
$$
Similarly, we obtain an orthogonal decomposition 
$$
l^2(G) = \mathbb{C} 1 \oplus \bigoplus_{n = 1}^\infty \; \bigoplus_{(i_1, \dots, i_n) \in I_n} l^2(G_{i_1} \setminus \{e\}) \otimes \cdots 
\otimes l^2(G_{i_n} \setminus \{e\})
$$
of the Hilbert space $ l^2(G) $. Let us write $ \mathbb{C}[G \setminus \{e\}]_{(j)} \subset \mathbb{C}[G\setminus \{e\}] $ for the 
subspace corresponding to tuples $ (i_1, \dots, i_n) \in I_n $ such that $ i_n \neq j $, and similarly write 
$ l^2(G \setminus \{e\})_{(j)} \subset l^2(G \setminus \{e\}) $ for the corresponding closed subspace. Note that we have 
$$ 
\mathbb{C}[G\setminus \{e\}] \cong \mathbb{C}[G\setminus \{e\}]_{(j)} \odot \mathbb{C}[G_j]
$$
and an analogous isomorphism for the Hilbert spaces. For $ j = 0, 1 $ consider the linear map 
$ T_j: \mathbb{C}[G \setminus \{e\}]_{(j)} \rightarrow \mathbb{C}[G/G_j] $ given by $ T_j(x) = x \odot 1 $ where we recall that 
$ \mathbb{C}[G/G_j] = \mathbb{C}[G] \odot_{\mathbb{C}[G_j]} \mathbb{C} $. Note that 
$$ 
\mathbb{C}[G/G_j] \cong \mathbb{C} 1 \oplus \bigoplus_{n = 1}^\infty \; \bigoplus_{(i_1, \dots, i_n) \in I^{(j)}_n} \mathbb{C}[G_{i_1}\setminus \{e\}] 
\odot \cdots \odot \mathbb{C}[G_{i_n}\setminus \{e\}]
$$ 
where 
$$
I^{(j)}_n = \{(i_1, \dots, i_n) \in \{0,1\}^n| i_k \neq i_{k + 1}\; \text{for all}\; k \;\text{and} \; i_n \neq j \}.
$$
The map $ T_j $ is a linear isomorphism onto the space $ \mathbb{C}[G/G_j \setminus \{e G_j\}] $. By the latter we mean the direct sum of 
all summands in the above decomposition of $ \mathbb{C}[G/H] $ except for the subspace $ \mathbb{C} 1 $ which represents the identity coset. 
Using analogous notation, this map extends to a unitary $ l^2(G)_{(j)} \rightarrow l^2(G/G_j \setminus \{e G_j \}) $ which we denote again 
by $ T_j $. \\
We define a bounded operator $ V_j^*: l^2(G) \rightarrow l^2(G/G_0) \oplus l^2(G/G_1) $ for $ j = 0,1 $ by $ V_j^* = T_k $ on 
$ l^2(G \setminus \{e\})_{(k)} $ for $ k = 0,1 $ and $ V_j^*(1) = 1 \in l^2(G/G_{1 - j}) $. That is, $ V_0^* $ differs from $ V_1^* $ 
only on the one-dimensional subspace $ \mathbb{C}1 \subset l^2(G) $ corresponding the trivial corepresentation. The operator 
$ V_j^* $ is an isometry with cokernel $ \im(V_j^*)^\bot = \mathbb{C}1 \subset l^2(G/G_j) $, and we define the Julg-Valette operator $ V_j $ 
to be the adjoint of $ V_j^* $. Geometrically, the operator $ V_j^* $ maps an edge to the endpoint furthest from the base point given by 
the identity coset in $ G/G_j $. \\
Now consider the odd operator $ \Gamma_j \in \LH(l^2(X)) $ given by  
$$
\Gamma_j = 
\begin{pmatrix} 
0 & V_j^* \\
V_j & 0 
\end{pmatrix} 
$$
with respect to the canonical decomposition $ l^2(X) = l^2(X^{(0)}) \oplus l^2(X^{(1)}) $. It is clear by construction that 
$ \Gamma_j \in \LH(l^2(X)) $ is self-adjoint and satisfies $ 1 - \Gamma_j^2 \in \KH(l^2(X)) $. \\
Let us write $ \lambda: l^2(X) \rightarrow M(C_0^\red(\DD(G)) \otimes l^2(X)) $ for the action of the Drinfeld double $ \DD(G) $ on $ l^2(X) $ 
and $ \ad_\lambda: \KH(l^2(X)) \rightarrow M(C^\red_0(\DD(G)) \otimes \KH(l^2(X)) $ for the associated adjoint coaction, compare 
\cite{NVpoincare}. The operator $ \Gamma_j $ almost commutes with the coaction $ \lambda $ on $ l^2(X) $ in the following sense. 
\begin{lemma} \label{almostequiv}
With the notation as above, we have 
$$
(C_0^\red(\DD(G)) \otimes 1)(1 \otimes \Gamma_j - \ad_\lambda(\Gamma_j)) \in C^\red_0(\DD(G)) \otimes \KH(l^2(X))
$$
for $ j = 0, 1 $. 
\end{lemma} 
\proof The coaction $ \lambda $ corresponds to coactions $ \alpha $ of $ : l^2(X) \rightarrow M(C_0(G) \otimes l^2(X)) $ and 
$ \nu: l^2(X) \rightarrow M(C^*_\red(G) \otimes l^2(X)) $ satisfying the Yetter-Drinfeld compatibility condition, see section \ref{secqg}. To 
check the assertion it suffices to show that $ \Gamma_j $ almost commutes with $ \alpha $ and $ \nu $ separately. The coaction $ \alpha $ 
corresponds to the canonical unital $ * $-representation of $ C^*_\max(G) = C^*_\max(G_0) * C^*_\max(G_1) $ on $ l^2(X) $. This in turn 
corresponds to a pair of unital $ * $-representations of $ C^*_\max(G_0) $ and $ C^*_\max(G_1) $ on $ l^2(X) $. It is easy to see that 
$ \Gamma_j $ commutes with the action of $ C^*_\max(G_j) $. Since $ \Gamma_0 $ and $ \Gamma_1 $ only differ by a finite rank operator, it 
follows that $ \Gamma_j $ commutes with the action of $ C^*_\max(G_{1 - j}) $ up to compact operators. This yields the claim for 
$ \alpha $. It is straightforward to check that $ \Gamma_j $ commutes strictly with the $ * $-representation of $ C_0(G) $ on $ l^2(X) $ 
induced by $ \nu $, and this yields the claim for $ \nu $. \qed \\
We conclude that the operator $ \Gamma_j $ together with the action of $ \mathbb{C} $ on $ l^2(X) $ by scalar multiplication defines a 
class $ \gamma \in KK^{\DD(G)}(\mathbb{C}, \mathbb{C}) $. Note that this class is independent of $ j $ since $ \Gamma_1 - \Gamma_0 $ is a 
finite rank operator. 
\begin{lemma} \label{gammaone}
We have $ \gamma = 1 $ in $ KK^{\DD(G)}(\mathbb{C}, \mathbb{C}) $. 
\end{lemma} 
\proof Let us use the operator $ \Gamma_0 $ to represent $ \gamma $. Then $ \gamma - 1 $ is represented by the graded Hilbert space 
$$ 
\H = \H_0 \oplus \H_1 =  l^2(X^{(0)}) \oplus (l^2(X^{(1)}) \oplus \mathbb{C}),  
$$
the action of $ \mathbb{C} $ by scalar multiplication and the self-adjoint odd unitary $ F_0 $ given as follows. We decompose $ \H $ into 
the direct sum $ \K \oplus \K^\bot $ where
$$ 
\K = \K_0 \oplus \K_1 = (\mathbb{C}1 \oplus \mathbb{C}1) \oplus (\mathbb{C}1 \oplus \mathbb{C}) 
\subset \H_0 \oplus \H_1
$$
and $ \mathbb{C}1 $ denotes the subspaces corresponding to the trivial corepresentation inside $ l^2(G/G_0), l^2(G/G_1) $ and $ l^2(G) $,
respectively. The operator $ F_0 $ preserves the decomposition $ \H = \K \oplus \K^\bot $, it agrees with $ \Gamma_0 $ on $ \K^\bot $, and 
on $ \K $ it interchanges $ 1 \in l^2(G/G_0) $ with $ 1 \in \mathbb{C} $ and $ 1 \in l^2(G/G_1) $ with $ 1 \in l^2(G) $. \\
As in lemma \ref{almostequiv} we study the $ \DD(G) $-action on $ \H $ in terms of the corresponding $ * $-representations of 
$ C^*_\max(G_0), C^*_\max(G_1) $ and $ C_0(G) $. It is obvious that $ F_0 $ commutes with the actions of $ C^*_\max(G_0) $ and $ C_0(G) $. 
Performing a symmetric construction with $ \Gamma_1 $ instead of $ \Gamma_0 $, we obtain an operator $ F_1 $ which commutes with the actions 
of $ C^*_\max(G_1) $ and $ C_0(G) $. Moreover $ F_0 = F_1 u $ where $ u $ is the self-adjoint unitary which is equal to the identity on 
$ \K^\bot $, interchanges $ 1 \in l^2(G/G_0) $ with $ 1 \in l^2(G/G_1) $ and fixes $ 1 \in l^2(G) $ and $ 1 \in \mathbb{C} $. If we conjugate 
the action of $ C^*_\max(G_1) $ by $ u $ and leave the action of $ C^*_\max(G_0) $ fixed, we thus obtain an action of $ C^*_\max(G) $ which 
commutes strictly with $ F_0 $. Choosing a path of unitaries in $ \LH(\K_0) $ connecting $ u $ and $ \id $ on $ \K_0 $ yields a homotopy 
between our cycle and a degenerate cycle. This shows $ \gamma - 1 = 0 $ as claimed. \qed \\
Our next aim is to define the dual-Dirac element $ \eta \in KK^{\DD(G)}(\mathbb{C}, \A_X) $. For a given $ \lambda \in E $ consider the 
operator
$$
\beta_\lambda(\mu) = \frac{c(\mu - \lambda)}{\sup(|\mu - \lambda|, 1/3)}
$$
in $ \LH(C_\tau(E)) = M(C_\tau(E)) $ where $ c $ denotes Clifford multiplication. By construction, the function $ 1 - \beta_\lambda^2 $ 
is supported on the ball $ B_{1/3}(\lambda) $ around $ \lambda $. The operator $ \beta_\lambda $ represents the Bott element in 
$ KK(\mathbb{C}, C_\tau(E)) $, we view it as an element $ [\beta_\lambda] $ in $ KK^{\DD(G)}(\mathbb{C}, C_\tau(E)) $ by considering the 
trivial action of $ \DD(G) $ on all ingredients of the corresponding cycle. Note that the class $ [\beta_\lambda] $ is independent of 
$ \lambda $. \\
Let us fix $ b = (1/2, 1/2) \in E $ and write $ \beta = \beta_b $. Since the action on $ C_\tau(E) $ is trivial, we obtain an element 
$ \tau_{C_\tau(E)}(\gamma) \in KK^{\DD(G)}(C_\tau(E), C_\tau(E)) $ by tensoring $ \gamma $ with $ C_\tau(E) $ in the standard way. The 
Kasparov product 
$$ 
\gamma \otimes_{\mathbb{C}} [\beta] = [\beta] \otimes_{C_\tau(E)} \tau_{C_\tau(E)}(\gamma) 
$$
is an element in $ KK^{\DD(G)}(\mathbb{C}, C_\tau(E)) $ which can be represented easily using the Julg-Valette operator $ \Gamma_0 $ from 
above. More precisely, the underlying Hilbert $ C_\tau(E) $-module is $ C_\tau(E) \cotimes l^2(X) $, the action of $ \mathbb{C} $ is by 
scalar multiplication, and the operator is given by 
$$
F_0 = \beta \cotimes 1 + ((1 - \beta^2)^{1/2} \cotimes 1)(1 \cotimes \Gamma_0), 
$$ 
see proposition 18.10.1 in \cite{Blackadar}. Using equivariant Morita invariance, we obtain a corresponding element in 
$ KK^{\DD(G)}(\mathbb{C}, C_\tau(E) \cotimes \hat{\KH}(l^2(X))) $ with underlying Hilbert module $ C_\tau(E) \cotimes \hat{\KH}(l^2(X)) $, 
the tautological left action of $ \mathbb{C} $, and the operator given by the same formula as $ F_0 $. \\
Let $ x_0, x_1 \in E $ be the points
$$
x_0 = (3/2, - 1/2), \qquad x_1 = (-1/2, 3/2),
$$
and let $ p_j \in M(C_0(G)) \subset \LH(l^2(G)) $ for $ j = 0,1 $ be the projection onto $ l^2(G \setminus \{e\})_{(j)} $. Note that 
$ p_0, p_1 \in M(C_0(G)) $ are central elements. For $ t \in [0,1] $ we define 
$ \beta^t \in \LH(C_\tau(E) \cotimes l^2(X)) = M(C_\tau(E) \cotimes \hat{\KH}(l^2(X))) $ by 
\begin{align*}
\beta^t &= \beta_{t x_0 + (1 - t) b} \cotimes (p_0 \oplus \id \oplus 0) \\
&\qquad + \beta_{t x_1 + (1 - t) b} \cotimes ((p_1 + (1 \otimes 1)) \oplus 0 \oplus \id)
\end{align*}
where we use the decomposition $ l^2(X) = l^2(G) \oplus l^2(G/G_0) \oplus l^2(G/G_1) $, and $ 1 \odot 1 \in \FH(l^2(G)) \subset \KH(l^2(G)) $ 
is the orthogonal projection onto $ \mathbb{C}1 \subset l^2(G) $. In addition, we define $ F_t \in \LH(C_\tau(E) \cotimes l^2(X)) $ for 
$ t \in [0,1] $ by 
$$ 
F_t = \beta^t + (1 - (\beta^t)^2)^{1/2}(1 \cotimes \Gamma_0), 
$$ 
for $ t = 0 $ this is of course compatible with our notation above. It is straightforward to check that $ (\beta^t)^2 $ commutes with 
$ \Gamma_0 $, and this implies $ F_t = F_t^* $. Moreover one verifies $ 1 - F_t^2 \in \KH(C_\tau(E) \cotimes l^2(X)) $ for all $ t \in [0,1] $ 
using that $ 1 - \beta_\lambda^2 $ is contained in $ C_\tau(E) $ for all $ \lambda \in E $ and $ 1 - \Gamma_0^2 \in \KH(l^2(X)) $. As in the 
proof of lemma \ref{almostequiv} we see that $ F_t $ is almost invariant under the action of $ \DD(G) $, that is 
$$ 
(C_0^\red(\DD(G)) \otimes 1)(1 \otimes F_t - \ad_\lambda(F_t)) \in C^\red_0(\DD(G)) \otimes \KH(C_\tau(E) \cotimes l^2(X))
$$
for all $ t \in [0,1] $, where we write again $ \lambda $ for the coaction on $ C_\tau(E) \otimes l^2(X) $ obtained from the coaction on 
$ l^2(X) $. \\
Let us now define the dual Dirac element $ \eta $. 
\begin{prop} \label{etadefinition}
The operator $ F_1 $ is contained in the multiplier algebra $ M(\A_X) $ inside $ \LH(C_\tau(E) \cotimes l^2(X)) $ and defines an element 
$ \eta \in KK^{\DD(G)}(\mathbb{C}, \A_X) $ in a canonical way. 
\end{prop} 
\proof Note that the canonical inclusion $ \A_X \rightarrow \LH(C_\tau(E) \cotimes l^2(X)) $ is nondegenerate, so that $ M(A_X) $ is indeed 
naturally a subalgebra of $ \LH(C_\tau(E) \cotimes l^2(X)) $. \\
To check that $ F_1 $ is a multiplier of $ \A_X $, it is convenient to write the operator in matrix form with respect to the decomposition 
$ l^2(X) = l^2(G/G_0) \oplus l^2(G/G_1) \oplus l^2(G) $. In a similar way one can represent the different linear subspaces of 
$ C_\tau(E) \cotimes \hat{\KH}(l^2(X)) $ in the definition of $ \A_X $. This reduces the argument to a number of straightforward verifications, 
the crucial point being that $ 1 - \beta_{x_j}^2 $ is supported on $ \Omega_j $ for $ j = 0, 1 $. Since this is completely analogous to the 
classical case treated in \cite{KSbuildings} we omit these verifications here. \\
We obtain the element $ \eta \in KK^{\DD(G)}(\mathbb{C}, \A_X) $ by considering $ \A_X $ as a Hilbert module over itself together with the 
operator $ F_1 \in \LH(\A_X) = M(\A_X) $ and the left action of $ \mathbb{C} $ by scalar multiplication. We have already seen above that 
$ F_1 $ is self-adjoint. To check the remaining properties of a Kasparov module, observe that 
$$ 
\KH(C_\tau(E) \cotimes l^2(X)) \cap M(\A_X) \subset \A_X = \KH(\A_X) 
$$ 
since $ \A_X $ acts nondegenerately on $ C_\tau(E) \cotimes l^2(X) $. Hence our previous considerations about the operators 
$ F_t \in \LH(C_\tau(E) \cotimes l^2(X)) $ show that $ 1 - F_1^2 \in \KH(\A_X) $ and that $ F_1 \in \LH(A_X) $ is almost invariant under the 
action of $ \DD(G) $. \qed \\
Let us compute the product $ \eta \otimes_{\A_X} D $. 
\begin{theorem} \label{gammaequalone}
Let $ D \in KK^{\DD(G)}(\A_X, \mathbb{C}) $ and $ \eta \in KK^{\DD(G)}(\mathbb{C}, \A_X) $ be the Dirac and dual Dirac elements defined above. 
Then 
$$
\eta \otimes_{\A_X} D = \id
$$
in $ KK^{\DD(G)}(\mathbb{C}, \mathbb{C}) $.
\end{theorem} 
\proof Let $ i: \A_X \rightarrow C_\tau(E) \cotimes \hat{\KH}(l^2(X)) $ be the canonical embedding and denote by 
$ [i] \in KK^{\DD(G)}(\A_X, C_\tau(E)) $ the corresponding class. Almost by construction, the operators 
$ F_t \in \LH(C_\tau(E) \cotimes l^2(X)) $ for $ t \in [0,1] $ define a homotopy between 
$ \eta \otimes_{\A_X} [i] $ and $ \gamma \otimes_\mathbb{C} [\beta] $. Hence if $ [\hat{\beta}] \in KK^{\DD(G)}(C_\tau(E), \mathbb{C}) $ denotes 
the inverse of the Bott element, then we obtain
$$
\eta \otimes_{\A_X} D = \eta \otimes_{\A_X} [i] \otimes_{C_\tau(E)} [\hat{\beta}]
= \gamma \otimes_\mathbb{C} [\beta] \otimes_{C_\tau(E)} [\hat{\beta}] = [\beta] \otimes_{C_\tau(E)} [\hat{\beta}] = \id 
$$
using lemma \ref{gammaone} and Bott periodicity. \qed

\section{The Baum-Connes conjecture and discrete quantum groups} \label{secbc}

In this section we first recall some elements of the categorical approach to the Baum-Connes conjecture developed by Meyer and Nest 
\cite{MNtriangulated}, \cite{MNhomalg1}, \cite{Meyerhomalg2}. In particular, we discuss the formulation of an analogue of the Baum-Connes 
conjecture for a certain class of discrete quantum groups proposed in \cite{Meyerhomalg2}. Using the results obtained in section 
\ref{secdirac} and \cite{Voigtbcfo} we will then show that free quantum groups satisfy this conjecture. \\
Let $ G $ be a discrete quantum group. The equivariant Kasparov category $ KK^G $ has as objects all separable $ G $-$ C^* $-algebras, 
and $ KK^G(A,B) $ as the set of morphisms between two objects $ A $ and $ B $. Composition of morphisms is given by the Kasparov product. 
The category $ KK^G $ is triangulated with translation automorphism $ \Sigma: KK^G \rightarrow KK^G $ given by the suspension 
$ \Sigma A = C_0(\mathbb{R}, A) $ of a $ G $-$ C^* $-algebra $ A $. Every $ G $-equivariant $ * $-homomorphism $ f: A \rightarrow B $ 
induces a diagram of the form
$$
\xymatrix{
\Sigma B  \;\; \ar@{->}[r] & C_f \ar@{->}[r] & A \ar@{->}[r]^f & B
}
$$
where $ C_f $ denotes the mapping cone of $ f $. Such diagrams are called mapping cone triangles. By definition, an exact triangle is a 
diagram in $ KK^G $ of the form $ \Sigma Q \rightarrow K \rightarrow E \rightarrow Q $ which is isomorphic to a mapping cone triangle. \\
Associated with the inclusion of the trivial quantum subgroup $ E \subset G $ we have the obvious restriction functor 
$ \res^G_E: KK^G \rightarrow KK^E = KK $ and an induction functor $ \ind_E^G: KK \rightarrow KK^G $. Explicitly, 
$ \ind_E^G(A) = C_0(G) \otimes A $ for $ A \in KK $ with coaction given by comultiplication on the copy of $ C_0(G) $. \\
We consider the following full subcategories of $ KK^G $, 
\begin{align*}
\TC_G &= \{A \in KK^G|\res^G_E(A) = 0 \in KK \} \\
\TI_G &= \{\ind_E^G(A)| A \in KK \}, 
\end{align*}
and refer to their objects as trivially contractible and trivially induced $ G $-$ C^* $-algebras, respectively. If there is no risk of 
confusion we will write $ \TC $ and $ \TI $ instead of $ \TC_G $ and $ \TI_G $. \\
The subcategory $ \TC $ is localising, and we denote by $ \bra \TI \ket $ the localising subcategory generated by $ \TI $. It follows from 
theorem 3.21 in \cite{Meyerhomalg2} that the pair of localising subcategories $ (\bra \TI \ket, \TC) $ in $ KK^G $ is complementary. That is, 
$ KK^G(P,N) = 0 $ for all $ P \in \bra \CI \ket $ and $ N \in \CC $, and every object $ A \in KK^G $ fits into an exact triangle 
$$
\xymatrix{
\Sigma N \; \ar@{->}[r] & \tilde{A} \ar@{->}[r] & A \ar@{->}[r] & N
}
$$
with $ \tilde{A} \in \bra \TI \ket $ and $ N \in \TC $. Such a triangle is uniquely determined up to isomorphism and depends functorially 
on $ A $. We will call the morphism $ \tilde{A} \rightarrow A $ a Dirac element for $ A $. \\
The localisation $ \mathbb{L}F $ of a homological functor $ F $ on $ KK^G $ at $ \TC $ is given by 
$$
\mathbb{L}F(A) = F(\tilde{A}) 
$$
where $ \tilde{A} \rightarrow A $ is a Dirac element for $ A $. By construction, there is an obvious map $ \mu_A: \mathbb{L}F(A) \rightarrow F(A) $ 
for every $ A \in KK^G $. \\
In the sequel we write $ G \ltimes_\max A $ and $ G \ltimes_\red A $ for the full and reduced crossed products of a $ G $-$ C^* $-algebra 
$ A $. 
\begin{definition} \label{defbc}
Let $ G $ be a discrete quantum group and consider the functor $ F(A) = K_*(G \ltimes_\red A) $ on $ KK^G $. 
We say that $ G $ satisfies the $ \TI $-Baum-Connes property with coefficients in $ A $ if the assembly map 
$$
\mu_A: \mathbb{L}F(A) \rightarrow F(A) 
$$
is an isomorphism. We say that $ G $ satisfies the $ \TI $-strong Baum-Connes property if $ \bra \TI \ket = KK^G $. 
\end{definition}
Clearly the $ \TI $-strong Baum-Connes property implies the $ \TI $-Baum-Connes property with coefficients in $ A $ for every 
$ G $-$ C^* $-algebra $ A $. We refer to \cite{MNtriangulated} for the comparison of definition \ref{defbc} with the usual 
formulation of the Baum-Connes conjecture in the case that $ G $ is a torsion-free discrete group. The $ \TI $-strong 
Baum-Connes property is equivalent to the assertion that $ G $ has a $ \gamma $-element and $ \gamma = 1 $ in this case. \\
In order to obtain the correct formulation of the Baum-Connes conjecture for an arbitrary discrete group $ G $ one has to replace 
the category $ \TI $ by the full subcategory $ \CI \subset KK^G $ of compactly induced $ G $-$ C^* $-algebras. Our results below show that for 
a free quantum group $ G $ the category $ \TI $ is already sufficient to generate $ KK^G $. It is natural to expect that all free 
quantum groups are torsion-free in the sense of \cite{Meyerhomalg2}. For free orthogonal quantum groups this is verified in 
\cite{Voigtbcfo}, but we have not checked the general case. \\
Let us briefly discuss some facts from \cite{MNhomalg1}, \cite{Meyerhomalg2} that will be needed in the sequel. 
We write $ \mathfrak{J} $ for the homological ideal in $ KK^G $ consisting of all $ f \in KK^G(A,B) $ such 
that $ \res^G_E(f) = 0 \in KK(A,B) $. By definition, $ \mathfrak{J} $ is the 
kernel of the exact functor $ \res^G_E: KK^G \rightarrow KK $. The ideal $ \mathfrak{J} $ is compatible with countable direct sums and has 
enough projective objects. The $ \mathfrak{J} $-projective objects in $ KK^G $ are precisely the retracts of compactly induced 
$ G $-$ C^* $-algebras. \\
A chain complex 
\begin{equation*}
\xymatrix{
\cdots \ar@{->}[r] & C_{n + 1} \ar@{->}[r]^{d_{n + 1}} & C_n \ar@{->}[r]^{d_n} & C_{n - 1} \ar@{->}[r] & \cdots 
}
\end{equation*}
in $ KK^G $ is $ \mathfrak{J} $-exact if 
\begin{equation*}
\xymatrix{
\cdots \ar@{->}[r] & KK(A, C_{n + 1}) \ar@{->}[r]^{\;\;\;(d_{n + 1})_*} & KK(A, C_n) \ar@{->}[r]^{(d_n)_*} & KK(A, C_{n - 1}) \ar@{->}[r] & \cdots 
}
\end{equation*}
is exact for every $ A \in KK $. \\
A $ \mathfrak{J} $-projective resolution of $ A \in KK^G $ is a chain complex 
\begin{equation*}
\xymatrix{
\cdots \ar@{->}[r] & P_{n + 1} \ar@{->}[r]^{d_{n + 1}} & P_n \ar@{->}[r]^{d_n} & P_{n - 1} \ar@{->}[r] & \cdots \ar@{->}[r]^{d_2} & P_1 \ar@{->}[r]^{d_1} & P_0
}
\end{equation*}
of $ \mathfrak{J} $-projective objects in $ KK^G $, augmented by a map $ P_0 \rightarrow A $ such that the augmented chain complex is 
$ \mathfrak{J} $-exact. \\
In the sequel we will denote by $ \hat{G} $ the dual compact quantum group of the discrete quantum group $ G $ determined 
by $ C^\red(\hat{G}) = C^*_\red(G)^\cop $. Note that this amounts to switching the comultiplication in the usual conventions. 
By definition, the representation ring of $ \hat{G} $ is the ring $ R(\hat{G}) = KK^{\hat{G}}(\mathbb{C}, \mathbb{C}) $ with multiplication 
given by Kasparov product. Note that $ R(\hat{G}) $ as an additive basis indexed by the irreducible corepresentations of $ G $. 
For every $ B \in KK^G $ there is a natural $ R(\hat{G}) $-module structure on 
$$ 
K(B) \cong K(\KH_G \otimes B) \cong KK^{\hat{G}}(\mathbb{C}, G \ltimes_\red B) 
$$ 
induced by Kasparov product. This module structure is important for doing homological algebra in $ KK^G $, see 
\cite{Meyerhomalg2} for further information. \\
Let us next record a basic result on induction and restriction for discrete quantum groups. 
\begin{prop} \label{indkk}
Let $ G $ be a discrete quantum group and let $ H \subset G $ be a quantum subgroup. Then 
$$
KK^G(\ind_H^G(A), B) \cong KK^H(A, \res^G_H(B))
$$
for every $ H $-$ C^* $-algebra $ A $ and every $ G $-$ C^* $-algebra $ B $. 
\end{prop}
\proof We describe the unit and counit of this adjunction, using the explicit description of the induced algebra $ \ind_H^G(A) $ from 
section \ref{secquantumsub}. We have an $ H $-equivariant $ * $-homomorphism $ \eta_A: A \rightarrow \ind_H^G(A) $ given by 
$ \eta_A(a) = \alpha(a) \in M(C_0(H) \otimes A) \subset M(C_0(G) \otimes A) $. This defines the unit 
$ \eta_A: A \rightarrow \res^G_H \ind_H^G(A) $ of the adjunction. \\
The counit $ \kappa_B: \ind_H^G \res^G_H(B) \cong C_0(G/H) \twisted B \rightarrow B $ of the adjunction is induced from the embedding 
$ C_0(G/H) \subset \KH(l^2(G/H)) $ followed by the $ \DD(G) $-equivariant Morita equivalence between $ \KH(l^2(G/H)) $ and $ \mathbb{C} $. 
One checks that 
$$
\xymatrix{
\ind_H^G(A) \ar@{->}[r]^{\!\!\!\!\!\!\!\!\!\!\!\!\!\!\! \ind(\eta_A)} & \; \ind_H^G \res^G_H \ind_H^G(A) 
\ar@{->}[r]^{\qquad \;\; \kappa_{\ind(A)}} & \;\ind_H^G(A)
}
$$
is equal to the identity in $ KK^G $ for any $ H $-$ C^* $-algebra $ A $, using that $ x \in \ind_H^G(A) $ is mapped to 
$ p_{[\epsilon]} \otimes x \in \KH(l^2(G/H)) \otimes \ind_H^G(A) $ under these maps. \\
Similarly the composition  
$$
\xymatrix{
\res^G_H(B) \ar@{->}[r]^{\!\!\!\!\!\!\!\!\!\!\!\!\!\!\! \eta(\res_B)} & \; \res^G_H \ind_H^G \res^G_H(B) 
\ar@{->}[r]^{\qquad \;\; \res_{\kappa(B)}} & \;\res^G_H(B)
}
$$
is the identity in $ KK^H $ for every $ G $-$ C^* $-algebra $ B $ since the minimal projection $ p_{[\epsilon]} \in \KH(l^2(G/H)) $ is 
$ H $-invariant and invariant under the adjoint coaction. This yields the claim. \qed \\
Finally, we need a description of the crossed product of homogeneous spaces for divisible quantum subgroups with respect to the conjugation 
coaction. Recall that the conjugation coaction $ \gamma: C_0(G) \rightarrow M(C^*_\red(G) \otimes C_0(G)) $ is given by 
$ \gamma(f) = \hat{W}^*(1 \otimes f) \hat{W} $ where $ \hat{W} = \Sigma W^* \Sigma \in M(C^*_\red(G) \otimes C_0(G)) $.
\begin{lemma} \label{homspacelemma} 
Let $ G $ be a discrete quantum group and let $ H \subset G $ be a divisible quantum subgroup. Then there exists a 
$ C_0(G)^\cop $-colinear $ * $-isomorphism 
$$
C_0(G)^\cop \ltimes_\red C_0(G/H) \cong C_0(G)^\cop \otimes C_0(G/H) 
$$
where the crossed product $ C_0(G)^\cop \ltimes_\red C_0(G/H) $ for the conjugation coaction is equipped with the dual coaction, and on 
$ C_0(G)^\cop \otimes C_0(G/H) $ we consider the comultiplication on the first tensor factor. 
\end{lemma} 
\proof Since we assume $ H $ to be divisible we find a central projection $ p \in M(C_0(G)) $ such that multiplication by $ p $ yields 
a $ * $-isomorphism from $ C_0(G/H) \subset C_b(G) $ to $ p C_0(G) $. Given an $ H $-equivariant $ * $-isomorphism 
$ C_0(G) \cong C_0(G/H) \otimes C_0(H) $, the projection $ p $ corresponds to $ 1 \otimes p_{[\epsilon]} $ where $ \epsilon \in \Irr(H) $ 
is the trivial corepresentation. Since $ p $ is central, the isomorphism $ C_0(G/H) \cong p C_0(G) $ commutes with the conjugation 
coaction. \\
We obtain 
\begin{align*} 
C_0(G)^\cop \ltimes_\red C_0(G/H) &\cong [(C_0(G) \otimes 1) \hat{W}^*(1 \otimes C_0(G)p) \hat{W}] \\
&\cong [\hat{W} (C_0(G) \otimes 1) \hat{W}^*(1 \otimes C_0(G) p)] \\
&\cong [\Sigma W^*(1 \otimes C_0(G)) W \Sigma(1 \otimes C_0(G)p)] \\
&= [\Delta^\cop(C_0(G))(1 \otimes C_0(G) p)] \\
&= [C_0(G) \otimes C_0(G) p] \\
&\cong C_0(G)^\cop \otimes C_0(G/H).
\end{align*}
It is straightforward to check that these identifications are compatible with the dual coaction on $ C_0(G)^\cop \ltimes_\red C_0(G/H) $ and 
comultiplication on the first tensor factor in $ C_0(G)^\cop \otimes C_0(G/H) $, respectively. \qed \\
Now let $ G_0 $ and $ G_1 $ be discrete quantum groups and let $ G = G_0 * G_1 $. In section \ref{secdirac} we constructed 
the Dirac element $ D \in KK^{\DD(G)}(\A_X, \mathbb{C}) $ using the $ C^* $-algebra $ \A_X $ of the quantum tree associated to $ G $. It was 
also shown that $ \A_X $ is $ KK^{\DD(G)} $-equivalent to a certain ungraded $ \DD(G) $-$ C^* $-algebra $ \P $. By slight abuse of notation, we 
will view $ D $ as an element in $ KK^{\DD(G)}(\P, \mathbb{C}) $ in the sequel. 
\begin{lemma} \label{Klemma} 
Let $ G_0, G_1 $ be discrete quantum groups and let $ G = G_0 * G_1 $ be their free product. The
Dirac element $ D \in KK(\P, \mathbb{C}) $ is invertible. 
\end{lemma} 
\proof By the definition of $ \P $ we have an extension 
$$
\xymatrix{
0 \ar@{->}[r] & I_0 \oplus I_1 \ar@{->}[r] & \P \ar@{->}[r] & \Sigma C(\Delta^1, C_0(G)) \ar@{->}[r] & 0 
}
$$
and using the $ KK $-equivalences between $ I_j $ and $ C_0(G/G_j) $ for $ j = 0,1 $ we obtain an exact sequence 
$$
\xymatrix{
 {0\;} \ar@{->}[r] \ar@{<-}[d] &
      K_1(\P) \ar@{->}[r] &
        K_0(C_0(G)) \ar@{->}[d]^\partial \\
   {0\;} \ar@{<-}[r] &
    {K_0(\P)}  \ar@{<-}[r] &
     {K_0(C_0(G/G_0)) \oplus K_0(C_0(G/G_1))} \\
}
$$
in $ K $-theory. Let us write $ e_r $ for the class in $ K_0(C_0(G)) $ corresponding to the irreducible corepresentation 
$ r \in \Irr(G) $. If $ r = \epsilon $ is the trivial corepresentation then we have 
$$
\partial(e_r) = (-e_{\tau_0(\epsilon)}, e_{\tau_1(\epsilon)})
$$
where $ e_r = p_\epsilon $ is the generator of $ K_0(C_0(G)) $ corresponding to $ \epsilon $, and 
$ e_{\tau_j(\epsilon)} $ is the generator corresponding to the class of $ \epsilon $ in $ K_0(C_0(G/G_j)) $. 
If $ r $ is a nontrivial corepresentation then we can write $ r = w t_j $ with $ w $ a possibly empty 
alternating word in $ \Irr(G_0) * \Irr(G_1) $ not ending in $ \Irr(G_j) $ and $ t_j \in \Irr(G_j) $ for $ j = 0 $ or $ j = 1 $. 
For $ j = 0 $ we find 
$$
\partial(e_r) = (- \dim(t_0) e_{\tau_0(w)}, e_{\tau_1(r)})
$$
where $ \tau_j(s) $ denotes the class corresponding $ s $ in the direct sum decomposition of $ C_0(G/G_j) $, and 
$ e_{\tau_j(s)} $ is the corresponding generator of $ K_0(C_0(G/G_j)) $. 
Similarly one obtains 
$$
\partial(e_r) = (-e_{\tau_0(r)}, \dim(t_1) e_{\tau_1(w)})
$$
for $ j = 1 $. In these formulas $ \dim(t_j) $ is the dimension of the Hilbert space underlying the corepresentation $ t_j $. \\
From this description it is easy to 
compute $ K_0(\P) = \mathbb{C} $ and $ K_1(\P) = 0 $. Using theorem \ref{gammaequalone} it follows that $ D $ induces an 
isomorphism in $ K $-theory. According to the 
universal coefficient theorem this yields the claim. \qed \\
The previous lemma can be strengthened as follows. 
\begin{lemma} \label{KDGlemma}
Let $ G = G_0 * G_1 $ be a free product of discrete quantum groups and set $ \hat{S} = C^*_\red(G) $. Then 
the Dirac element $ D \in KK^{\hat{S}}(\P, \mathbb{C}) $ is invertible. 
\end{lemma} 
\proof Note that the coaction of $ C^*_\red(G) $ is obtained by restriction from the coaction of the Drinfeld double $ C_0^\red(\DD(G)) $. 
If we write $ S = C_0(G)^\cop $, then by Baaj-Skandalis duality the assertion is equivalent to saying that 
$ S \ltimes_\red D \in KK^S(S \ltimes_\red \P, S \ltimes_\red \mathbb{C}) $ is invertible. According to theorem \ref{gammaequalone} 
it suffices in fact to show that the map 
$ KK^S(S \ltimes_\red \P, S \ltimes_\red \P) \rightarrow KK^S(S \ltimes_\red \P, S \ltimes_\red \mathbb{C}) $ 
induced by $ D $ is an isomorphism. We have an extension
$$
\xymatrix{
0 \ar@{->}[r] & S \ltimes_\red (I_0 \oplus I_1) \ar@{->}[r] & S \ltimes_\red \P \ar@{->}[r] & 
S \ltimes_\red \Sigma C(\Delta^1, C_0(G)) \ar@{->}[r] & 0 
}
$$
of $ S $-$ C^* $-algebras with equivariant completely positive splitting, and inserting this extension into the first variable yields a 
$ 6 $-term exact sequence in $ KK^S $. According to lemma \ref{homspacelemma} we see that $ S \ltimes_\red I_j $ and 
$ S \ltimes_\red \Sigma C(\Delta^1, C_0(G)) $ are $ S $-equivariantly isomorphic to algebras of the form $ S \otimes B $ where 
$ B \in KK $ is in the bootstrap class. Due to proposition \ref{indkk} it therefore suffices to show that 
$ S \ltimes_\red D \in KK(S \ltimes_\red \P, S \ltimes_\red \mathbb{C}) $ induces an isomorphism in $ K $-theory. \\
Using lemma \ref{homspacelemma} we obtain a commutative diagram 
\begin{equation*}
\xymatrix{
S \ltimes_\red (I_0 \oplus I_1) \ar@{->}[r] \ar@{->}[d]^\cong & S \ltimes_\red \P \ar@{->}[r] & S \ltimes_\red \Sigma C(\Delta^1, C_0(G)) 
\ar@{->}[r] \ar@{->}[d]^\cong & \Sigma (S \ltimes_\red (I_0 \oplus I_1)) \ar@{->}[d]^\cong \\ 
S \otimes (I_0 \oplus I_1) \ar@{->}[r] & S \otimes \P \ar@{->}[r] & S \otimes \Sigma C(\Delta^1, C_0(G)) \ar@{->}[r] & 
\Sigma (S \otimes (I_0 \oplus I_1))
}
\end{equation*}
where the rows are exact triangles. Since the category $ KK^S $ is triangulated there exists an invertible morphism 
$ S \ltimes_\red \P \rightarrow S \otimes \P $ in $ KK^S $ completing the above commutative diagram. 
In particular we obtain a diagram
\begin{equation*}
\xymatrix{
K_0(S \ltimes_\red \P) \ar@{->}[r]^{S \ltimes_\red D} \ar@{->}[d]^{\cong} 
& K_0(S \ltimes_\red \mathbb{C}) \ar@{->}[d]^{\cong} \\ 
K_0(S \otimes \P) \ar@{->}[r]^{S \otimes D} & K_0(S \otimes \mathbb{C})
}
\end{equation*}
where the lower horizontal arrow is an isomorphism according to lemma \ref{Klemma}. To check commutativity of this diagram 
recall that we have a natural $ R(\hat{G}) $-module structure on the $ K $-groups $ K(S \ltimes_\red B) $ for any $ B $, where 
$ R(\hat{G}) \cong KK^{\hat{S}}(\mathbb{C}, \mathbb{C}) $ is the representation ring of $ \hat{G} $. 
A generator of the $ R(\hat{G}) $-module $ K_0(S \otimes \P) $ is represented by the projection 
$ p_\epsilon \otimes i(p_{[\epsilon]},0) $ where $ p_\epsilon \in S $ and $ p_{[\epsilon]} \in C_c(G/G_0) $ are the projections 
corresponding to the trivial corepresentation and $ i: C_0(G/G_0) \oplus C_0(G/G_1) \rightarrow I_0 \oplus I_1 \rightarrow \P $ is the 
canonical morphism. Since $ p_{[\epsilon]} $ is invariant under the conjugation action of $ C_c(G) $ one finds that 
$ p_\epsilon \ltimes i(p_{[\epsilon]},0) $ is a generator of the $ R(\hat{G}) $-module $ K_0(S \ltimes_\red \P) $. 
Under $ S \ltimes_\red D $ this element is mapped to $ p_\epsilon \in K_0(S) = K_0(S \ltimes_\red \mathbb{C}) $. 
Since all the maps in the above diagram are $ R(\hat{G}) $-linear this finishes the proof. \qed \\
We shall now show that the $ \TI $-strong Baum-Connes property is inherited by free products, compare \cite{Oyonotrees}, \cite{Tutrees}. 
\begin{theorem} \label{bcfreeprod}
Let $ G_0 $ and $ G_1 $ be discrete quantum groups satisfying the $ \TI $-strong Baum-Connes property and let $ G = G_0 * G_1 $ be 
their free product. Then $ G $ satisfies the $ \TI $-strong Baum-Connes property as well. 
\end{theorem}
\proof Let $ A \in KK^G $ be an arbitrary $ G $-algebra. Using theorem 3.6 in \cite{NVpoincare} we see that the braided tensor product 
$ \P \twisted A $ is an extension of $ G $-$ C^* $-algebras induced from $ G_0, G_1 $ and the trivial quantum subgroup. Therefore the 
$ \TI $-strong Baum-Connes property for $ G_0 $ and $ G_1 $ implies that $ \P \twisted A $ is contained in $ \bra \TI_G \ket $. \\
It remains to prove that $ D \twisted A \in KK^G(\P \twisted A, A) $ is an isomorphism. Due to theorem \ref{gammaequalone} 
it is enough to show that the map $ KK^G(\P \twisted A, \P \twisted A) \rightarrow KK^G(\P \twisted A, A) $ induced by $ D \twisted A $ 
is an isomorphism. \\
By construction of $ \P $ we have an extension 
$$
\xymatrix{
0 \ar@{->}[r] & I_0 \twisted A \oplus I_1 \twisted A \ar@{->}[r] & \P \twisted A \ar@{->}[r] & \Sigma C(\Delta^1, C_0(G)) \twisted A \ar@{->}[r] & 0 
}
$$
of $ G $-$ C^* $-algebras with equivariant completely positive splitting.  
Using the six-term exact sequence in $ KK $-theory and proposition \ref{indkk} we may thus reduce 
the problem to showing that the maps $ KK^{G_j}(B, \P \twisted A) \rightarrow KK^{G_j}(B, A) $ 
and $ KK(B, \P \twisted A) \rightarrow KK(B, A) $ are isomorphisms for all $ B $ in $ KK^{G_j} $ and $ KK $, respectively. 
Since $ G_j $ satisfies the $ \TI $-strong Baum-Connes property, every $ B \in KK^{G_j} $ is contained in $ \bra \TI_{G_j} \ket $. 
Hence it is in fact enough to show that $ KK(B, \P \twisted A) \rightarrow KK(B, A) $ 
is an isomorphism for all $ B \in KK $, or equivalently that $ D \twisted A \in KK(\P \twisted A, A) $ 
is invertible. However, this follows from lemma \ref{KDGlemma} and the functorial properties of the braided tensor product, 
compare \cite{NVpoincare}. \qed 
\begin{lemma} \label{bcsubgroup}
Let $ G $ be a discrete quantum group satisfying the $ \TI $-strong Baum-Connes property. 
If $ H \subset G $ is a divisible quantum subgroup, then $ H $ satisfies the $ \TI $-strong Baum-Connes property
as well. 
\end{lemma}
\proof Let $ A \in KK^H $ be given. Since $ G $ satisfies the $ \TI $-strong Baum-Connes property we have 
$ \ind_H^G(A) \in \bra \TI_G \ket $ in $ KK^G $. The category $ \TI_G $ is generated by $ G $-$ C^* $-algebras of the form 
$ C_0(G) \otimes B $ for $ B \in KK $ with coaction given by the comultiplication of $ C_0(G) $. Since by assumption $ H \subset G $ is 
divisible we have $ C_0(G) = C_0(H) \otimes C_0(H \backslash G) $ as $ H $-$ C^* $-algebras. This implies 
$ \res^G_H(\bra \TI_G \ket) \subset \bra \TI_H \ket $ and in particular $ \res^G_H \ind_H^G(A) \in \bra \TI_H \ket $ in $ KK^H $. Using 
the explicit description of induced $ C^* $-algebras in section \ref{secquantumsub} we see that the $ H $-$ C^* $-algebra $ A $ is a 
retract of $ \res^G_H \ind_H^G(A) $. That is, there are morphisms $ A \rightarrow \res^G_H \ind_H^G(A) $ and $ \res^G_H \ind_H^G(A) \rightarrow A $  which compose to the identity on $ A $. Since $ \bra \TI_H \ket $ is closed under retracts this shows $ A \in \bra \TI_H \ket $ and 
yields the claim. \qed \\
We are now ready to prove the main results of this paper. 
\begin{theorem} \label{bcfu}
Let $ n > 1 $ and $ Q \in GL_n(\mathbb{C}) $. Then the free unitary quantum group $ \mathbb{F}U(Q) $ 
satisfies the $ \TI $-strong Baum-Connes property. 
\end{theorem}
\proof We consider first the case $ Q \in GL_2(\mathbb{C}) $, and without loss of generality we may assume that $ Q $ is positive. Then 
$ Q $ can be written in the form
$$
Q = r
\begin{pmatrix}
q & 0 \\
0 & q^{-1}
\end{pmatrix}
$$ 
for some positive real number $ r $ and $ q \in (0,1] $, and $ \mathbb{F}U(Q) $ is isomorphic to $ \mathbb{F}U(P) $ where 
$$
P = 
\begin{pmatrix}
0 & q^{1/2} \\
-q^{-1/2} & 0
\end{pmatrix}
$$
satisfies $ P \overline{P} = -1 $. We remark that the free orthogonal quantum group $ \mathbb{F}O(P) $ is isomorphic 
to the dual of $ SU_q(2) $. \\
It is well-known that $ \mathbb{Z} $ satisfies the $ \TI $-strong Baum-Connes property \cite{HKatmenable}, and according to \cite{Voigtbcfo} 
the same is true for $ \mathbb{F}O(P) $. Hence due to proposition \ref{fuinc} the free product $ \mathbb{F}O(P) * \mathbb{Z} $ satisfies the 
$ \TI $-strong Baum-Connes property as well. Using proposition \ref{fuinc} and proposition \ref{bcsubgroup} we conclude that $ \mathbb{F}U(P) $
satisfies the $ \TI $-strong Baum-Connes property. This yields the assertion for $ Q \in GL_2(\mathbb{C}) $. \\
Now let $ Q \in GL_n(\mathbb{C}) $ be arbitrary. According to corollary 6.3 in \cite{BdRV} the dual of $ \mathbb{F}U(Q) $ is monoidally 
equivalent to the dual of $ \mathbb{F}U(R) $ for a suitable matrix $ R \in GL_2(\mathbb{C}) $. Hence the claim follows from the invariance of 
the $ \TI $-strong Baum-Connes property under monoidal equivalence, see theorem 8.6 in \cite{Voigtbcfo}. \qed \\
Combining this with theorem \ref{bcfreeprod} and the results in \cite{Voigtbcfo} we obtain the following theorem. 
\begin{theorem} \label{bcfree}
Let $ G $ be a free quantum group of the form 
$$
G = \mathbb{F}U(P_1) * \cdots * \mathbb{F}U(P_k) * \mathbb{F}O(Q_1) * \cdots * \mathbb{F}O(Q_l)
$$
for matrices $ P_j \in GL_{m_i}(\mathbb{C}) $ with $ m_i > 1 $ for all $ i $ and matrices $ Q_j \in GL_{n_j}(\mathbb{C}) $ with $ n_j > 2 $ 
such that $ Q_j \overline{Q_j} = \pm 1 $ for all $ j $. Then $ G $ satisfies the $ \TI $-strong Baum-Connes property. 
\end{theorem}

\section{$ K $-theory for free quantum groups} \label{secktheory}

In this section we discuss the main applications of our results. We shall establish an analogue of the Pimsner-Voiculescu exact 
sequence for free quantum groups and compute their $ K $-theory. In addition we discuss some consequences concerning idempotents 
in reduced $ C^* $-algebras of free quantum groups and the $ \gamma $-element described in \cite{Vergniouxtrees}. \\
Let $ G $ be a discrete quantum group and $ A $ a $ G $-$ C^* $-algebra. For $ s \in \Irr(G) $ we define $ s_* \in KK(A,A) $ 
as the composition 
$$
\xymatrix{
A \ar@{->}[r]^{\!\!\!\!\!\!\!\!\!\!\!\!\!\!\!\!\alpha} & \; C(G)_{\overline{s}} \otimes A
\ar@{->}[r]^{\quad \; \; \; \simeq} & A 
}
$$
where the first arrow is obtained by composing the coaction $ \alpha: A \rightarrow M(C_0(G) \otimes A) $ with the projection onto the matrix 
block corresponding to $ \overline{s} $ inside $ C_0(G) $. The second arrow is induced by the Morita equivalence 
$ C(G)_{\overline{s}} \cong \LH(\H_{\overline{s}}) \sim_M \mathbb{C} $. In this way we obtain a ring homomorphism 
$ R(\hat{G}) \rightarrow KK(A, A) $, where we recall that $ R(\hat{G}) = KK^{\hat{G}}(\mathbb{C}, \mathbb{C}) $ denotes the representation ring 
of the dual compact quantum group $ \hat{G} $. In particular, the representation ring acts on the $ K $-theory of $ A $, and this 
action agrees in fact with the action defined in section \ref{secbc}. \\
Using this notation we shall now formulate the Pimsner-Voiculescu exact sequence for free quantum groups. 
\begin{theorem} \label{pv} 
Let $ G $ be a free quantum group of the form 
$$
G = \mathbb{F}U(P_1) * \cdots * \mathbb{F}U(P_k) * \mathbb{F}O(Q_1) * \cdots * \mathbb{F}O(Q_l)
$$
for matrices $ P_j \in GL_{m_i}(\mathbb{C}) $ with $ m_i > 1 $ for all $ i $ and matrices
$ Q_j \in GL_{n_j}(\mathbb{C}) $ with $ n_j > 2 $ such that $ Q_j \overline{Q_j} = \pm 1 $ for all $ j $. \\
Then $ G $ is $ K $-amenable, and therefore the natural map 
$$
K_*(G \ltimes_\max A) \rightarrow K_*(G \ltimes_\red A)
$$
is an isomorphism for every $ G $-$ C^* $-algebra $ A $. Moreover there is an exact sequence 
$$
\xymatrix{
 \bigoplus_{j = 1}^{k + 2l} K_0(A) \ar@{->}[r]^{\quad \sigma} \ar@{<-}[d] &
      K_0(A) \ar@{->}[r]^{\!\!\!\!\!\! \iota_*} &
        K_0(G \ltimes_\max A) \ar@{->}[d] \\
   K_1(G \ltimes_\max A) \ar@{<-}[r]^{\quad \;\; \iota_*} &
    K_1(A) \ar@{<-}[r]^{\!\!\! \sigma} &
    \bigoplus_{j = 1}^{k + 2l} K_1(A) \\
}
$$
for the $ K $-theory of the crossed product. Here $ \sigma $ is the map 
$$
\sigma = \bigoplus_{i = 1}^k ((u_i)_* - m_i) \oplus ((\overline{u}_i)_* - m_i) \oplus \bigoplus_{j = 1}^l (v_j)_* - n_j 
$$
where we write $ u_i $ and $ v_j $ for the fundamental corepresentations of $ \mathbb{F}U(P_i) $ and $ \mathbb{F}O(Q_j) $, 
respectively. In addition $ \iota: A \rightarrow G \ltimes_\max A $ denotes the canonical inclusion. 
\end{theorem}
\proof Since full and reduced crossed products agree for trivially induced actions, theorem \ref{bcfree} implies that 
$ G \ltimes_\max A \rightarrow G \ltimes_\red A $ induces an isomorphism in $ KK $ for every $ G $-$ C^* $-algebra $ A $, compare 
\cite{Voigtbcfo}. This means that $ G $ is $ K $-amenable. \\
Consider the homological ideal $ \mathfrak{J} $ in $ KK^G $ given by the kernel of the restriction 
functor $ \res^G_E: KK^G \rightarrow KK $. We shall now construct a $ \mathfrak{J} $-projective resolution of length $ 1 $ 
for every $ G $-$ C^* $-algebra $ A $. \\
Note first that the representation ring $ R(\hat{G}) = KK^{\hat{G}}(\mathbb{C}, \mathbb{C}) $ can be identified with the free
algebra $ \mathbb{Z}\bra u_1, \overline{u_1}, \dots, u_k, \overline{u_k}, v_1, \dots, v_l \ket $ generated by the fundamental corepresentations 
$ u_i $ of $ \mathbb{F}U(P_i) $, their conjugates, and the fundamental corepresentations $ v_j $ of $ \mathbb{F}O(Q_j) $. 
We may also view $ R(\hat{G}) $ as the tensor algebra $ TV $ over 
$$ 
V = \bigoplus_{i = 1}^k (\mathbb{Z} u_i \oplus \mathbb{Z} \overline{u}_i) \oplus  \bigoplus_{j = 1}^l \mathbb{Z} v_j. 
$$
Note that for $ B = C_0(G) $ the $ R(\hat{G}) $-module structure on $ K(B) \cong R(\hat{G}) $ is given by multiplication, and for 
$ B = \mathbb{C} $ this module structure is induced from the homomorphism $ \epsilon: TV \rightarrow \mathbb{Z} $ given by 
$ \epsilon(u_i) = m_i = \epsilon(\overline{u}) $ and $ \epsilon(v_j) = n_j $. \\
For $ r \in \Irr(G) $ we define an element $ T_r \in KK^{\DD(G)}(C_0(G), C_0(G)) $ by 
$$
\xymatrix{
C_0(G) \ar@{->}[r]^{\!\!\!\!\!\!\!\!\!\!\!\!\!\!\!\!\Delta} & \; C_0(G) \otimes C(G)_r
\ar@{->}[r]^{\quad \; \; \; \simeq} & C_0(G) 
}
$$
where the first arrow is given by composition of the comultiplication of $ C_0(G) $ with the projection onto the matrix block 
$ C(G)_r \cong \KH(\H_r) $ corresponding to $ r $. The second morphism is given by the canonical $ \DD(G) $-equivariant Morita equivalence 
between $ C_0(G) \otimes \KH(\H_r) \cong \KH(C_0(G) \otimes \H_r) $ and $ C_0(G) $. It is straightforward to check that the induced map 
$ T_r: K_0(C_0(G)) \rightarrow K_0(C_0(G)) $ identifies with right multiplication with $ \overline{r} $ under the identification 
$ K_0(C_0(G)) \cong R(\hat{G}) $. \\
Let us consider the diagram
$$
\xymatrix{
0 \ar@{->}[r] & P_1 \ar@{->}[r]^{\delta} & P_0 \ar@{->}[r]^{\lambda} & \mathbb{C} \ar@{->}[r] & 0
}
$$
in $ KK^{\DD(G)} $ where 
$$
P_1 = \bigoplus_{j = 1}^{2k + l} C_0(G), \qquad P_0 = C_0(G) 
$$
and the morphisms are defined as follows. The arrow $ \lambda: C_0(G) \rightarrow \KH_G \simeq \mathbb{C} $ is given by the regular 
representation. Moreover 
$$ 
\delta = \bigoplus_{i = 1}^k (T_{u_i} - m_i \id) \oplus (T_{\overline{u_i}} - m_i \id) \oplus \bigoplus_{j = 1}^l (T_{v_j} - n_j \id)
$$ 
where $ T_r $ for $ r \in \Irr(G) $ are the morphisms defined above. \\
By applying $ K $-theory to the above diagram we obtain the sequence of $ TV $-modules 
$$
\xymatrix{
0 \ar@{->}[r] & TV \otimes V \ar@{->}[r]^{\;\;d} & TV \ar@{->}[r]^{\epsilon} & \mathbb{Z} \ar@{->}[r] & 0 
}
$$
where $ \epsilon $ is the augmentation homomorphism from above, and $ d $ is given by 
$$
d(1 \otimes u_i) = \overline{u_i} - m_i, \qquad d(1 \otimes \overline{u_i}) = u_i - m_i, \qquad d(1 \otimes v_j) = v_j - n_j 
$$
on the canonical $ TV $-basis of $ TV \otimes V $. It is easy to write down a $ \mathbb{Z} $-linear contracting 
homotopy for this complex. Moreover this homotopy can be lifted to $ KK^{\hat{S}} $ where $ \hat{S} = C^*_\red(G) $ since 
the conjugation coaction is diagonal with respect to the canonical direct sum decomposition of $ C_0(G) $, 
and each matrix block $ C(G)_r $ for $ r \in \Irr(G) $ is $ \hat{S} $-equivariantly Morita equivalent to $ \mathbb{C} $. \\
Now let $ A $ be a $ G $-$ C^* $-algebra. Then taking braided tensor products yields a sequence
$$
\xymatrix{
0 \ar@{->}[r] & P_1 \twisted A \ar@{->}[r]^{\!\!\!\delta \twisted \id} & P_0 \twisted A 
\ar@{->}[r]^{\!\!\lambda \twisted \id} & \mathbb{C} \twisted A \ar@{->}[r] & 0
}
$$
in $ KK^G $ which is split exact in $ KK $. Taking into account that $ C_0(G) \twisted A \cong C_0(G) \otimes A $ is $ \mathfrak{J} $-projective 
and $ \mathbb{C} \twisted A \cong A $, we conclude that the above sequence defines a $ \mathfrak{J} $-projective resolution of $ A $. 
We have thus shown that every object of $ KK^G $ has a projective resolution of length $ 1 $. Moreover $ KK^G(A, B) = 0 $ for 
every $ A \in KK^G $ and all $ \mathfrak{J} $-contractible objects $ B $. Indeed, this relation holds for $ A \in \TI $ due to 
proposition \ref{indkk}, and hence for all objects in $ KK^G = \bra \TI \ket $ according to theorem \ref{bcfree}. \\
Consider the homological functor $ F $ defined on $ KK^G $ with values in abelian groups given by $ F(A) = K(G \ltimes_\max A) $. 
According to theorem 4.4 in \cite{MNhomalg1} we obtain a short exact sequence 
$$
\xymatrix{
0 \ar@{->}[r] & \mathbb{L}_0 F_*(A) \ar@{->}[r] & K_*(G \ltimes_\max A) \ar@{->}[r] & \mathbb{L}_1 F_{* - 1}(A) \ar@{->}[r] & 0 
}
$$
for all $ n $. Moreover we have 
\begin{align*}
\mathbb{L}_0 F_*(A) &= \coker(F_*(\delta): K_*(A) \rightarrow K_*(A)) \\
\mathbb{L}_1 F_*(A) &= \ker(F_*(\delta): K_*(A) \rightarrow K_*(A)).
\end{align*}
It is easy to check that $ F_*(\delta) $ identifies with the map $ \sigma $. We can thus splice the above short exact sequences together to 
obtain the desired six-term exact sequence. Finally, it is straightforward to identify the map $ K_*(A) \rightarrow K_*(G \ltimes_\max A) $
in this sequence with the map induced by the inclusion. \qed \\
Using theorem \ref{pv} we may compute the $ K $-groups of free quantum groups. Let us first explicitly state the result in the unitary case. 
\begin{theorem} \label{futheorem}
Let $ n > 1 $ and $ Q \in GL_n(\mathbb{C}) $. Then the natural homomorphism $ C^*_\max(\mathbb{F}U(Q)) \rightarrow C^*_\red(\mathbb{F}U(Q)) $ 
induces an isomorphism in $ K $-theory and 
$$
K_0(C^*_\max(\mathbb{F}U(Q))) = \mathbb{Z}, \qquad K_1(C^*_\max(\mathbb{F}U(Q))) = \mathbb{Z} \oplus \mathbb{Z}.
$$
These groups are generated by the class of $ 1 $ in the even case and the classes of $ u $ and $ \overline{u} $ in the odd case. 
\end{theorem}
\proof Let us abbreviate $ G = \mathbb{F}U(Q) $. According to theorem \ref{pv} the natural map $ K_*(C^*_\max(G)) \rightarrow K_*(C^*_\red(G)) $ 
is an isomorphism. To compute these groups we consider the projective resolution constructed in the proof of theorem \ref{pv}. The resulting 
long exact sequence in $ K $-theory takes the form 
$$
\xymatrix{
 {\mathbb{Z}\;} \ar@{->}[r] \ar@{<-}[d] &
      K_0(C^*_\max(G)) \ar@{->}[r] &
        0 \ar@{->}[d] \\
   {\mathbb{Z}^2\;} \ar@{<-}[r] &
    {K_1(C^*_\max(G))}  \ar@{<-}[r] &
     {0} \\
}
$$
which yields 
$$
K_0(C^*_\max(G)) = \mathbb{Z}, \qquad K_1(C^*_\max(G)) = \mathbb{Z} \oplus \mathbb{Z}.
$$
Indeed, the Hopf-$ C^* $-algebra $ C^*_\max(G) $ has a counit, and therefore the upper left horizontal map in the above diagram 
is an isomorphism. We conclude in particular that the even $ K $-group $ K_0(C^*_\max(G)) $ is generated by the class of $ 1 $. \\
It remains to determine the generators of $ K_1(C^*_\max(G)) = K_1(C^*_\max(\mathbb{F}U(Q))) $. 
Assume first $ Q \overline{Q} = \pm 1 $ and consider the automorphism $ \tau $ of $ C^*_\max(G) $ determined by 
$ \tau(u) = Q \overline{u} Q^{-1} $. Note that 
$$
\tau(Q \overline{u} Q^{-1}) = Q \tau(\overline{u}) Q^{-1} = Q \overline{Q} u \overline{Q}^{-1} Q^{-1} = u 
$$
so that $ \tau $ is a well-defined $ * $-homomorphism such that $ \tau^2 = \id $. 
Moreover $ \tau $ is compatible with the comultiplication. \\
Using $ \tau $ we obtain for any $ G $-$ C^* $-algebra $ A $ a new $ G $-$ C^*$-algebra $ A^{\hat{\tau}} $ by 
considering $ A = A^{\hat{\tau}} $ with the coaction $ \alpha^{\hat{\tau}}(x) = (\hat{\tau} \otimes \id)\alpha(x) $, 
where $ \hat{\tau} $ is the automorphism of $ C_0(G) $ dual to $ \tau $. This automorphism can be viewed as an  
equivariant $ * $-homomorphism $ \hat{\tau}: C_0(G) \rightarrow C_0(G)^{\hat{\tau}} $.  
On the level of $ K $-theory we obtain a commutative diagram 
\begin{equation*}
\xymatrix{
0 \ar@{->}[r] & K_1(C^*_\red(G)) \ar@{->}[r] \ar@{->}[d]^{\tau_*} & \mathbb{Z} \oplus \mathbb{Z}
\ar@{->}[r] \ar@{->}[d]^\sigma & \mathbb{Z} \ar@{->}[r] \ar@{->}[d]^\id & K_0(C^*_\red(G)) \ar@{->}[r] \ar@{->}[d]^{\tau_*} & 0 \\ 
0 \ar@{->}[r] & K_1(C^*_\red(G)) \ar@{->}[r] & \mathbb{Z} \oplus \mathbb{Z} \ar@{->}[r] & \mathbb{Z} \ar@{->}[r] & K_0(C^*_\red(G)) \ar@{->}[r] & 0
}
\end{equation*}
where $ \sigma(k,l) = (l,k) $ is the flip map. It follows in particular that the elements in $ \mathbb{Z} \oplus \mathbb{Z} $ corresponding to 
$ u $ and $ \overline{u} $ are obtained from each other by applying $ \sigma $, that is, 
$$
[u] = (p,q), \qquad [\overline{u}] = (q,p)
$$
for some $ p, q \in \mathbb{Z} $. \\
Note that $ C^*_\max(\mathbb{F}U(Q)) = C^*_\max(\mathbb{F}U(|Q|)) $ where 
$ |Q| = (Q^* Q)^{1/2} $ is the absolute value of $ Q $. Using this identification 
we obtain a $ * $-homomorphism $ \rho: C^*_\max(\mathbb{F}U(Q)) \rightarrow C^*_\max(\mathbb{Z}) $ by choosing a suitable corner 
of the generating matrix. More precisely, if $ u $ is the generating matrix of $ C^*_\max(\mathbb{F}U(|Q|)) $ then 
we may set $ \rho(u_{11}) = z $ and $ \rho(u_{ij}) = \delta_{ij} $ else. On the level of $ K_1 $ 
the map $ \rho $ sends $ [u] $ to $ 1 \in \mathbb{Z} = K_1(C^*_\max(\mathbb{Z})) $ and $ [\overline{u}] $ to $ -1 $. 
Under the canonical projection $ C^*_\max(\mathbb{F}U(Q)) \rightarrow C^*_\max(\mathbb{F}O(Q)) $ the 
classes $ [u] $ and $ [\overline{u}] $ are both mapped to $ 1 \in \mathbb{Z} = K_1(C^*_\max(\mathbb{F}O(Q))) $. \\
Consider the $ * $-homomorphism $ \rho \oplus \pi: C^*_\max(\mathbb{F}U(Q)) \rightarrow C^*_\max(\mathbb{Z}) \oplus C^*_\max(\mathbb{F}O(Q)) $. 
The map induced by $ \rho \oplus \pi $ on the level of $ K_1 $ may be written as a matrix 
$$
A = \begin{pmatrix}
a & b \\
c & d 
\end{pmatrix} 
\in M_2(\mathbb{Z}).
$$
If we set 
$$
U = 
\begin{pmatrix}
p & q \\
q & p
\end{pmatrix}, 
\qquad 
M = 
\begin{pmatrix}
1 & -1 \\
1 & 1 
\end{pmatrix} 
$$
then our above argument shows
$$
A U = M. 
$$
In particular $ \det(U) = p^2 - q^2 = (p + q)(p - q) $ divides $ 2 = \det(M) $. We conclude that 
$ \det(U) = \pm 1 $ or $ \det(U) = \pm 2 $. The latter 
is impossible since either $ p + q $ and $ p - q $ are both divisible by two, or none of the factors is. Hence $ \det(U) = \pm 1 $, and 
this implies that $ u $ and $ \overline{u} $ generate $ K_1(C^*_\max(\mathbb{F}U(Q)) $. This yields the claim 
for $ Q $ satisfying $ Q \overline{Q} = \pm 1 $, and hence our assertion holds for all quantum groups $ \mathbb{F}U(Q) $ with 
$ Q \in GL_2(\mathbb{C}) $. \\
Now let $ Q \in GL_n(\mathbb{C}) $ be an arbitrary matrix. Without loss of generality we may assume that $ Q $ is positive. 
Then we find $ P \in GL_{2n}(\mathbb{C}) $ satisfying $ P \overline{P} \in \mathbb{R} $ such that $ Q $ is contained as a matrix block 
inside $ P $, and an associated surjective $ * $-homomorphism $ \alpha: C^*_\max(P) \rightarrow C^*_\max(Q) $. Similarly we find a 
$ * $-homomorphism $ \beta: C^*_\max(\mathbb{F}U(Q)) \rightarrow C^*_\max(\mathbb{F}U(R)) $ where $ R \in M_2(\mathbb{C}) $ is a matrix 
block contained in $ Q $. \\
By our above results we know that $ K_1(C^*_\max(\mathbb{F}U(P)) $ and $ K_1(C^*_\max(\mathbb{F}U(R))) $ are generated by 
the fundamental unitaries and their conjugates. Moreover, on the level of $ K_1 $ these generators are preserved
under the composition $ \beta_* \alpha_* $. Hence $ \beta_* \alpha_* = \id $, and we conclude that both $ \alpha_* $ and 
$ \beta_* $ are invertible. It follows that $ K_1(C^*_\max(Q)) $ is generated by the fundamental unitary and its conjugate as well. \qed \\
Combining theorem \ref{futheorem} with the results in \cite{Voigtbcfo} yields the following theorem. 
\begin{theorem} \label{freetheorem}
Let $ G $ be a free quantum group of the form 
$$
G = \mathbb{F}U(P_1) * \cdots * \mathbb{F}U(P_k) * \mathbb{F}O(Q_1) * \cdots * \mathbb{F}O(Q_l)
$$
for matrices $ P_i \in GL_{m_i}(\mathbb{C}) $ with $ m_i > 1 $ for all $ i $ and $ Q_j \in GL_{n_j}(\mathbb{C}) $ with $ n_j > 2 $ for all 
$ j $ such that $ Q_j \overline{Q_j} = \pm 1 $. \\ 
Then the $ K $-theory $ K_*(C^*_\max(G)) = K_*(C^*_\red(G)) $ of $ G $ is given by 
$$
K_0(C^*_\max(G)) = \mathbb{Z}, \qquad K_1(C^*_\max(G)) = \mathbb{Z}^{2k} \oplus \mathbb{Z}^l, 
$$
and these groups are generated by the class of $ 1 $ in the even case and the canonical unitaries in the odd case. 
\end{theorem}
\proof We may either proceed as in theorem \ref{futheorem}, or apply the $ K $-theory exact sequence for free products, using induction on 
the number of factors in $ G $. The case of a single factor is treated in theorem \ref{futheorem} and \cite{Voigtbcfo}, respectively. For the 
induction step we observe
$$ 
K_0(C^*_\max(H) * C^*_\max(K)) = (K_0(C^*_\max(H)) \oplus K_0(C^*_\max(K))/\bra ([1], -[1]) \ket 
$$
and 
$$
K_1(C^*_\max(H) * C^*_\max(K)) = K_1(C^*_\max(H)) \oplus K_1(C^*_\max(K))
$$
for arbitrary discrete quantum groups $ H $ and $ K $, see \cite{Cuntzfreeproduct}. The claim on the generators follows 
from these identifications and the above mentioned results. \qed \\
We shall write $ \mathbb{F}U(n) = \mathbb{F}U(1_n) $ and $ \mathbb{F}O(n) = \mathbb{F}O(1_n) $ if $ 1_n \in GL_n(\mathbb{C}) $ is the identity 
matrix. As a consequence of theorem \ref{freetheorem} we obtain the following result concerning idempotents in the reduced group $ C^* $-algebras 
of free quantum groups. The proof is the same as in the classical case of free groups, see also \cite{Voigtbcfo}. 
\begin{cor}
Let $ G $ be a free quantum group of the form 
$$
G = \mathbb{F}U(m_1) * \cdots * \mathbb{F}U(m_k) * \mathbb{F}O(n_1) * \cdots * \mathbb{F}O(n_l)
$$
for some integers $ m_i, n_j $ with $ n_j > 2 $ for all $ j $. Then $ C^*_\red(G) $ does not contain nontrivial idempotents. 
\end{cor}
In \cite{Vergniouxtrees} an analogue of the Julg-Valette element for free quantum groups is constructed. Our results imply that this element 
is homotopic to the identity, see the argument in \cite{Voigtbcfo}. 
\begin{cor} 
Let $ G $ be a free quantum group of the form 
$$
G = \mathbb{F}U(P_1) * \cdots * \mathbb{F}U(P_k) * \mathbb{F}O(Q_1) * \cdots * \mathbb{F}O(Q_l)
$$
for matrices $ P_i \in GL_{m_i}(\mathbb{C}) $ and $ Q_j \in GL_{n_j}(\mathbb{C}) $ with $ n_j > 2 $ for all $ j $ such that 
$ Q_j \overline{Q_j} = \pm 1 $. Then the Julg-Valette element for $ G $ is equal to $ 1 $ in  $ KK^G(\mathbb{C}, \mathbb{C}) $. 
\end{cor}

\bibliographystyle{plain}

\bibliography{cvoigt}

\end{document}